\documentclass[a4paper,12pt]{article}

\usepackage{amsfonts}
\usepackage{amsmath}
\usepackage{amssymb}
\usepackage{algorithm}
\usepackage{booktabs}
\usepackage{graphicx}
\usepackage{color}
\usepackage{enumerate}

\newtheorem{theorem}{Theorem}

\newtheorem{corollary}{Corollary}

\newtheorem{remark}{Remark}
\newtheorem{problem}{Problem}

\def\vol{\mathrm{vol}}

\def\trace{\mathrm{trace}}
\def\one{\mathbb{I}}

\usepackage{amsmath}
\usepackage{amssymb}
\usepackage{amsfonts}
\usepackage{graphicx}

\newcommand{\K}{\mathcal{K}}
\newcommand{\B}{\mathcal{B}}

\newcommand{\U}{\mathcal{U}}
\newcommand{\V}{\mathcal{V}}

\newcommand{\Poly}[1]{P_{#1}}
\newcommand{\unif}[1]{ {\mathbb{U}}_{#1} }

\newcommand{\SoS}[1]{\Sigma_{#1}}

\newcommand{\ped}[1]{_{\rm #1}}

\usepackage{ulem}

\usepackage[figuresright]{rotating}

\title{Simple Approximations of Semialgebraic Sets and their Applications to Control}

\author{Fabrizio Dabbene$^1$ \and Didier Henrion$^{2,3,4}$ \and Constantino Lagoa$^5$}

\footnotetext[1]{CNR-IEIIT; c/o Politecnico di Torino; C.so Duca degli Abruzzi 24, Torino; Italy}
\footnotetext[2]{CNRS; LAAS; 7 avenue du colonel Roche, F-31400 Toulouse; France.}
\footnotetext[3]{Universit\'e de Toulouse; LAAS; F-31400 Toulouse; France.}
\footnotetext[4]{Faculty of Electrical Engineering, Czech Technical University in Prague, Technick\'a 2, CZ-16626 Prague, Czech Republic.}
\footnotetext[5]{Electrical Engineering Department, The Pennsylvania State University, University Park, PA 16802, USA. }

\date{\today}

\begin{document}

\maketitle

\begin{abstract}
Many uncertainty sets encountered in control systems analysis and design can be expressed in terms of semialgebraic sets, that is as the intersection of sets described by means of polynomial inequalities. Important examples are for instance the solution set of linear matrix inequalities or the Schur/Hurwitz stability domains. 
These sets often have very complicated shapes (non-convex, and even non-connected), which renders very difficult their manipulation.
It is therefore of considerable importance to find simple-enough approximations of these sets,  able to capture their main characteristics while maintaining a low level of complexity.
For these reasons, in the past years several convex approximations, based for instance on hyperrectangles, polytopes, or ellipsoids have been proposed. 

In this work, we move a step further, and propose possibly non-convex approximations, based 
on a small volume polynomial superlevel set  of a single positive polynomial of given degree. 
We show how these sets can be easily approximated by minimizing the $L^{1}$ norm of the polynomial over the semialgebraic set, subject to positivity constraints. Intuitively, this corresponds to the trace minimization heuristic  
commonly encounter in minimum volume ellipsoid problems. 
From a computational viewpoint, we design a hierarchy of  linear matrix inequality 
problems to generate these approximations, and we  provide theoretically rigorous convergence results, in the sense that the hierarchy of outer approximations converges in volume (or, equivalently, almost  everywhere and almost uniformly) to the original set.

Two main applications of the proposed approach are considered. The first one aims at reconstruction/approximation of sets from a finite number of  samples. In the second one, we show how the concept of polynomial superlevel set
can be used to generate samples uniformly distributed on a given semialgebraic set.
The efficiency of the proposed approach is demonstrated by different numerical examples.
\end{abstract}


\newpage
\section{Introduction}
\label{sec:intro}

In this paper, we address the problem of how to determine ``simple'' approximations of semialgebraic sets
{in Euclidean space}, and we show how these approximations can be exploited  to address several problems of interest in systems and control. To be more precise, given a
 set
\begin{equation}
\label{K-set}
\K\doteq \{x \in {\mathbb R}^n : g_i(x) \geq 0, \:i=1,2,\ldots,m\}
\end{equation}
which is compact, with non-empty interior and
described by {given real multivariate polynomials $g_i(x), i=1,2,\ldots,m$}, {and a compact set $\B\supset \K$},
we aim at determining a  so-called {\it polynomial superlevel set} {(PSS)} 
\begin{equation}\label{pss}
\U(p) \doteq \{ x \in \B : p(x) \geq 1 \}.
\end{equation}
that constitutes a good {outer} approximation of the set $\K$ of interest and converges {strongly} to $\K$
when increasing the degree of the real multivariate polynomial $p$ to be found.

In particular, the proposed PSS is based on an easily computable polynomial approximation of the indicator function of the set $\K$. In the paper, we show that suitable approximations of the indicator function can be obtained by solving a convex optimization problem whose constraints are linear matrix inequalities (LMIs) and that, as the degree of the approximation increases, one converges in $L^1$-norm, almost uniformly and almost everywhere to the indicator function of the semialgebraic set $\K$ of interest.
Moreover, the set approximations provided in this paper can be thought as a direct generalization of classical ellipsoidal set approximations, in the sense that if  second degree approximations are used,  we exactly recover well-known approaches.

The main motivation for the problem addressed in the paper is the fact that semialgebraic sets are frequently encountered in control. As an example, consider the Hurwitz or Schur stability regions of a polynomial. It is a well-known fact that the these regions are semialgebraic sets in the coefficient space.  The polynomial inequalities that define these stability sets can be derived from
well-known {algebraic stability criteria.}
Another classical example of semialgebraic sets arising in control are LMI
feasibility sets, {also called spectrahedra}. Indeed, LMI sets are 
(convex) basic semialgebraic sets. To see this, consider the LMI set
\[
	\K\ped{LMI}\doteq\{x \in {\mathbb R}^n : F(x)=F_{0}+F_{1}x_{1}+\cdots +F_{n}x_{n}\succeq 0\}
\]
where the matrix $F(x)$ has size $m\times m$, {and observe that
a vector $x$ belongs to $\K\ped{LMI}$ if and only if all the coefficients of the univariate
polynomial
\[
s \mapsto \mathrm{det}\:(sI_m+F(x)) =  g_1(x)+g_2(x)s+\cdots+g_m(x)s^{m-1}+s^m
\]
are nonnegative, i.e. $x$ belongs to the set $\K$ is defined in (\ref{K-set}), where the polynomials $g_i(x)$
are by construction sums of principal minors of the matrix $F(x)$.} 
The approach taken in this paper is the following: given the set $\K$, we search for a minimum volume PSS
that contains the set $\K$. Since there is {in general} no analytic formula for the volume of a semialgebraic set,
in terms of the coefficients of the polynomials defining the set\footnote{{See however reference \cite{MorSha:10}
which explains how explicit formulas can be obtained with discriminants in exceptional cases.}},
 it is {very challenging} to solve this optimization problem {locally, let alone} globally.
Instead, the main contribution of this paper is to  describe and justify analytically and geometrically a computationally tractable heuristic based on \mbox{$L^1$-norm} or trace
minimization. Second, we show that the same approach can be employed to obtain the largest (in terms of the $L^1$ surrogate for the volume) PSS inscribed in $\K$.
Moreover, it is shown how the ideas put forth in this paper can be used to address two important problems: i) reconstruction/approximation of a (possibly non-semialgebraic) set from samples belonging to it, and ii) uniform
generation of samples distributed over a semialgebraic set. Examples of applications in a systems analysis and controller design context are also provided. 

The work presented in this paper is an extension of the preliminary results in the conference papers~\cite{DabHen:13} and~\cite{DaHeLa:14}, and it provides a more in depth analysis of both theoretical and implementation aspects.
In particular, with respect to ~\cite{DabHen:13}, the present manuscript contains more detailed proofs of the theoretical results, provides detailed algorithmic descriptions, and introduces inner PSS approximations.
Similarly, the results on random sample generation of~\cite{DaHeLa:14} are here described in more details, and an algorithm is provided. Finally, all examples in the paper are new, and more control oriented applications are considered.

\subsection{Previous work and related literature}

The idea of approximating overly complicated sets by introducing simpler and easy manageable geometrical shapes is surely not new, it has a very long history, and it arises
in different research fields such as optimization, system identification and control.
In particular, in the systems and control community, the most common approach is to introduce  \textit{outer bounding sets}, that is sets of minimum size which are guaranteed to contain the set to be approximated.
For instance, in the context of robust filtering, set-theoretic state estimators for uncertain nonlinear dynamic systems have been proposed in \cite{AlBrCa:05,DuWaPo:01,ElGCal:01,ShaTu:97}. 
These strategies adopt a set-membership approach \cite{GaTeVi:99,Schweppe:73}, and construct (the smallest) compact set guaranteed to bound the system states that are consistent with the measured output and the norm-bounded uncertainty. The most common geometrical shape adopted in these work is the ellipsoidal one, for the double reason that it has a very simple description -- the center and the shape matrix are sufficient to provide a complete characterization -- and that its determination usually can be formulated as a convex (usually quadratic) optimization problem. The use of ellipsoidal sets in the state estimation problems was introduced in the pioneering work \cite{Schweppe:73} and used by many different authors from then on; see, for example, \cite{DuWaPo:01,ElGCal:01}. 
Outer approximation also arise in the context of robust fault detection problems
(e.g., see \cite{IBPAG:09}) and of reachability analysis of nonlinear and/or hybrid systems \cite{HwStTo:03,KurVar:00}.
Similarly, \textit{inner} approximations  are employed in nonlinear programming \cite{NesNem:94}, in the solution of design centering problems \cite{WojVla:93} and for fixed-order controller design \cite{HenLou:12}.
In this case, one aims at constructing the set largest size inscribed in the set of interest.

Besides ellipsoids, other shapes have been considered in the recent literature.
The use of polyhedrons was proposed in \cite{KunLyc:85} to obtain an increased estimation accuracy, while zonotopes have been also recently studied in \cite{AlBrCa:05,GuNgZa:03}.
{In \cite{CePiRe:11} a heuristic based on polynomial
optimization and convex relaxations is proposed for computing small volume polytopic outer approximations of a compact semialgebraic set.}
More recent works, like for instance \cite{BeViLa:07,HenLou:12,MaLaBo:05}, employ sets defined by semialgebraic conditions.
The closest approach to the one proposed in our paper can be found in \cite{MaLaBo:05}, in which the authors use polynomial sum-of-squares (SOS) programming to address the problem of fitting given data with a convex polynomial, seen as a natural extension of quadratic polynomials and ellipsoids. Convexity of the polynomial is ensured by enforcing that 
its Hessian is matrix SOS, and volume minimization is indirectly
enforced by increasing the curvature of the polynomial.
In \cite{BeViLa:07} the authors propose moment-based relaxations for the separation
and covering problems with semialgebraic sets, thereby also extending the
classical ellipsoidal sets used in data fitting problems.

Recently, the authors of \cite{DaLaSh:10} have proposed an approach  based on randomization, which constructs convex approximations of generic nonconvex sets which are \textit{neither inner nor outer}, but they enjoy some specific probabilistic properties. In this context, an approximation is considered to be reliable if it contains ``most'' of the points in the given set
with prescribed high probability. The key tool in this framework is the generation of random samples inside
the given set, and the construction of a convex set containing these samples.

\subsection{The sequel}

The paper is organized as follows. In Section~\ref{sec:prob_stat} the notation used in this paper is introduced and the central problem addressed in this paper is defined. In order to be able to numerically solve the set approximation problem of interest, in Section~\ref{sec:l1_form} a related polynomial optimization problem is introduced and numerical methods for solving it are described in Section~\ref{sec:sos}.  In Section~\ref{sec:inner}, we discuss how the results in this paper can be used to find inner approximations of semialgebraic sets. A first set of numerical examples is provided in Section~\ref{sec:examples1}. Using the central results on set approximation mentioned above, in Section~\ref{sec:finite_samples} we address the  problem of reconstructing a set from a finite number of points in its interior. In Section~\ref{sec:sampling} we provide algorithms for uniform sample generation in semialgebraic sets and in Section~\ref{sec:conclusion} some closing remarks are provided.

\section{Problem statement} \label{sec:prob_stat}

Before a description of the main problem addressed is provided, we introduce the basic notation that is used throughout the paper.

\subsection{Notation}

The notation $A \succ 0$ ($\succeq 0$) means that the symmetric matrix $A$ is positive definite (semidefinite), and given two matrices $A$ and $B$ we write $A \succeq B$ whenever $A-B \succeq 0$. Given a set $\K\subset {\mathbb R}^n$,
its indicator function is defined as
\begin{equation}
\label{indic}
\one_\K(x) \doteq \begin{cases}
1 & \text{ if } x \in \K \\ 
0 & \text{ if } x \notin \K 
\end{cases}
\end{equation}
and its volume or, more precisely, the Lebesgue measure of $\K$, is denoted by
\[
\vol\:\K \doteq \int_{\K} dx = \int_{{\mathbb R}^n} \one_\K(x) dx.
\]
The set of all real coefficient polynomials of {degree} less than or equal to $d$ is denoted by $\Poly{d}$. The monomial basis for this set is represented by the (column) vector $\pi_d \in \Poly{d}$, so that any $p \in  \Poly{d}$ can be expressed in the following form
\[
\label{eq:Gram}
p(x) = \pi_d^T(x) p = \pi^T_{\lceil d/2 \rceil}(x) P \pi_{\lceil d/2 \rceil}(x)
\]
where $p$ is a real (column) vector\footnote{Note that we use $p$ to denote both the polynomial
and the vector of its coefficients whenever no ambiguity is possible.} and $P$ is a symmetric matrix of appropriate size,
often referred to as Gram matrix. Also, we denote by  $\SoS{2d}$ the set of polynomials $p \in \Poly{2d}$
that can be represented as sums of squares of other polynomials, i.e.
\[
p = \sum\limits_{k=1}^{n_p} p_k^2, \quad p_k\in \Poly{d}, \:k=1,\ldots,n_p.
\]
Finally, given a polynomial $p \in \Poly{d}$, define its  $L^1$ norm over a compact set $\B$, denoted by $L^1_\B$ or just $L^1$ when the set~$\B$ used is clear from the context, as
\[
\|p\|_1 \doteq \int_\B p(x) dx.
\]

\subsection{Problem Statement}

With the notation defined above, we are now ready to define the central problem in this paper. We consider the
 basic semialgebraic set $\K$ defined in (\ref{K-set}), which is assumed to be compact and with a non-empty interior.

As discussed in the Introduction, the set $\K$ has typically a complex description in terms of its defining polynomials
(e.g. coming from physical measurements and/or estimations). For this reason, we aim at finding a ``simpler'' approximation
of this set which has enough degrees of freedom to capture its characteristics. {This approximation is
the polynomial superlevel set (PSS) $\U(p)$ defined in (\ref{pss}) in terms of a real multivariate polynomial
$p \in \Poly{d}$ of given degree $d$. This degree controls the complexity of the approximation.}
Among the family of possible PSS that can be constructed, we 
search for the one that provides the set $\U(p)$ of minimum volume while containing the set of interest $\mathcal{K}$, hence capturing most of its the geometric features. 
Formally, we define the following optimization problem
\begin{problem}[Minimum volume outer {PSS}]
\label{prob1}
Given {$d \in {\mathbb N}$ and} a compact semialgebraic set $\K$, 
find a polynomial $p \in \Poly{d}$ whose PSS $\U(p)$ is 
of minimum volume and  contains $\K$. That is, solve the following optimization problem
\begin{equation}\label{opt}
\begin{array}{rcll}
{v^*_d} & {\doteq} & \displaystyle \inf_{p \in \Poly{d}} & \vol\:\U(p) \\
&& \mathrm{s.t.} & \K \subseteq \U(p). \\
\end{array}
\end{equation}
\end{problem}


Note that this problem can be viewed as the natural extension of the problem of computing the minimum volume ellipsoid containing $\K$. Indeed, if $\K$ is convex and the polynomial $p$ is quadratic ($d=2$), then the infimum of problem \eqref{opt} is attained, and the optimal set  $\U(p)$ is given by the unique (convex) ellipsoid of minimum volume that
contains $\K$, called L\"owner-John ellipsoid.
In particular, if $\K$ is the convex-hull of a finite set of points, this ellipsoid can be computed by convex optimization, see e.g. \cite[\S 4.9]{BenNem:01}.

We remark however that, for $d$ greater than $2$,  the optimization problem \eqref{opt} is nonlinear and semi-infinite, in the sense that the optimization is over the finite-dimensional vector space $\Poly{d}$,
but subject to an infinite number of constraints, necessary to cope with the set inclusion.
\begin{theorem}\label{inftomin}
{The sequence of infima of problem (\ref{opt}) monotically converges from above to $\vol\:\K$, i.e. for all $d\geq 1$
it holds $v^*_d \geq v^*_{d+1}$ and $\lim_{d\to\infty} v^*_d = \vol\:\K$.}
\end{theorem}

\noindent
\textbf{Proof:} 
{As in \cite[Section 3.2]{HeLaSa:09}, let $x \mapsto d(x,\K)$ be the Euclidean distance to set $\K$ and with $\epsilon_k > 0$ let $\K_\epsilon := \{x \in \B \: : \: d(x,\K) < \epsilon_k\}$ be an open bounded outer approximation of $\K$, so that $\B \backslash \K_k$ is closed with $\lim_{k\to\infty} \epsilon_k = 0$. By Urysohn's Lemma \cite[Section 12.1]{RoyFit:10}  there is a sequence of continuous functions $(f_k)_{k \in \mathbb N}$ with $f_k : \B \to [0,1]$ such that $f_k = 0$ on $\B \backslash \K_k$ and $f_k = 1$ on $\K$.
In particular, notice that $\vol\:\K \leq \vol\:\U(f_k) \leq \vol\:\K + \vol\:\K_k \backslash \K$ and since $\lim_{k\to\infty}
\vol\:\K_k \backslash \K = 0$ it holds $\lim_{k\to\infty} \vol\:\U(f_k) = \vol\:\K$.

By the Stone-Weierstrass Theorem \cite[Section 12.3]{RoyFit:10}
we can approximate $f_k$ uniformly on $\B$ by a sequence of polynomials $(p'_{k,d})_{d\in\mathbb N}$ with $p'_{k,d} \in \Poly{d}$,
i.e. $\sup_{x \in \B} |f_k(x)-p'_{k,d}(x)| < \epsilon'_d$ with $\lim_{d\to\infty} \epsilon'_d = 0$. Defining $p_{k,d}:=p'_{k,d}+2\epsilon'_d$,
the sequence of polynomials $(p_{k,d})_{d\in \mathbb N}$ converges uniformly to $f_k$ from above, i.e.
$p_{k,d} \geq f_k$ on $\B$ and $\lim_{d\to\infty} \sup_{x \in \B} |f_k(x)-p_{k,d}(x)|=0$.
This implies that $v^*_{k,d}:=\vol\:\U(p_{k,d}) \geq \vol\:\U(f_k)$ and $\lim_{d\to\infty} v^*_{k,d} = \vol\:\U(f_k)$.
Recalling $\lim_{k\to\infty} \vol\:\U(f_k) = \vol\:\K$, it follows that $\lim_{k,d\to\infty} v^*_{k,d}  = \vol\:\K$
which proves, up to extracting a subsequence indexed by $d$, the existence of a minimizing sequence of polynomials for optimization problem (\ref{opt}).

Finally, the inequality $v^*_d \geq v^*_{d+1}$ readily follows from the inclusion $\Poly{d} \subset \Poly{d+1}$.
}
$\Box$.

Before introducing the approach we propose for the solution of Problem \ref{prob1}, in the next subjection we briefly recall some recent results which are closely related to the problem considered in this paper, for the special case of homogeneous polynomials.

\subsection{Remark on a convex conic formulation}

In this section, we summarize existing results for the case when the polynomial $q\in\Poly{d}$ defined as $q(x)\doteq 2-p(x)$ is assumed to be a \textit{homogeneous} polynomial, or form, of
even degree $d=2\delta$ in $n$ variables. 
First note that, with this change of notation\footnote{ %
The polynomial $q\doteq 2-p$ is introduced {because} the results in  \cite{Lasserre:15}
are derived for  sublevel sets, {not superlevel sets}.}
 the {PSS} $\U(p)$ corresponds to the {unit} sublevel set ${\V(q)\doteq}\{ x \in {\mathbb R}^n : q(x) \leq 1 \}$ of the polynomial $q$.
In \cite[Lemma 2.4]{Lasserre:15}
it is proved that, when $q$ is homogeneous, the volume function
\[
q \mapsto \vol\:\V(q)
\]
is convex in $q$. The proof of this statement relies on the striking observation \cite{MorSha:10} that
\[
\vol\:\V(q) = C_d \int_{{\mathbb R}^n} e^{-q(x)}dx
\]
where $C_d$ is a constant depending only on $d$.
Note also that boundedness of $\V(q)$ implies that $q$ is nonnegative,
since if there is a point $x_0 \in {\mathbb R}^n$ such that $q(x_0)<0$, and hence $x_0 \in \V(q)$,
then by homogeneity of $q$ it follows that $q(\lambda x_0)=\lambda^{2\delta} q(x_0)<0$ for all $\lambda$
and hence $\lambda x_0 \in \V(p)$ for all $\lambda$ which contradicts
boundedness of $\V(p)$. This implies that problem (\ref{opt}), once restricted to nonnegative forms,
is a convex optimization problem.

Moreover, in \cite[Lemma 2.4]{Lasserre:15} explicit expressions are given for the
first and second order derivatives of the volume function, in terms
of the moments
\begin{equation}\label{mom}
\int_{{\mathbb R}^n} x^{\alpha} e^{-q(x)}dx
\end{equation}
for $\alpha \in {\mathbb N}^n$, $|\alpha|\leq 2d$.
In an iterative algorithm solving convex problem (\ref{opt}),
one should then be able to compute repeatedly
and quickly integrals of this kind, arguably a difficult task.
Moreover, when $q$ is not homogeneous, we do not know under which
conditions on $q$ the function $\vol\:\V(q)$ is convex in $q$.

Motivated by these considerations, in the remainder of this paper we propose a simpler approach
to the solution problem (\ref{opt}), which is not restricted to
forms, and which does not require the potentially intricate
numerical computation of moments (\ref{mom}) {of exponentials of homogeneous polynomials}.
The introduction of this approach is motivated by its {analogy} with the well-known trace heuristic for ellipsoidal approximation, and it consists of approximating the volume by means of the $L^1$-norm of the polynomial $p$.

\section{$L^1$-norm minimization} \label{sec:l1_form}

It is assumed that a ``simple set'' $\B \subset {\mathbb R}^n$ containing $\K$ is known. 
{By ``simple'' we mean that analytic expressions of the moments of the Lebesgue measure on $\B$ should be
available, so that integration of polynomials can be carried out readily.}
In the following, we  assume that the set $\B$ is an $n$-dimensional  hyperrectangle of the form
 \begin{equation}
 \label{Bab} 
{\B = [a,b]} \doteq\{x \in {\mathbb R}^n : a_i\leq x_i \leq b_i, \: i=1,2,\ldots,n \}
 \end{equation}
{with $a$ and $b$ given vectors of ${\mathbb R}^n$.}
This is a very mild assumption since, given a semialgebraic set like the set $\K$ above, one can easily compute an hyperrectangle containing it; see Section~\ref{sec:Bab} for details. We note that more complex sets $\B\supseteq\K$ can be considered, provided that integration of polynomials over it is easily done. 

\begin{figure}[!ht]
\centerline{
\includegraphics[width=8cm]{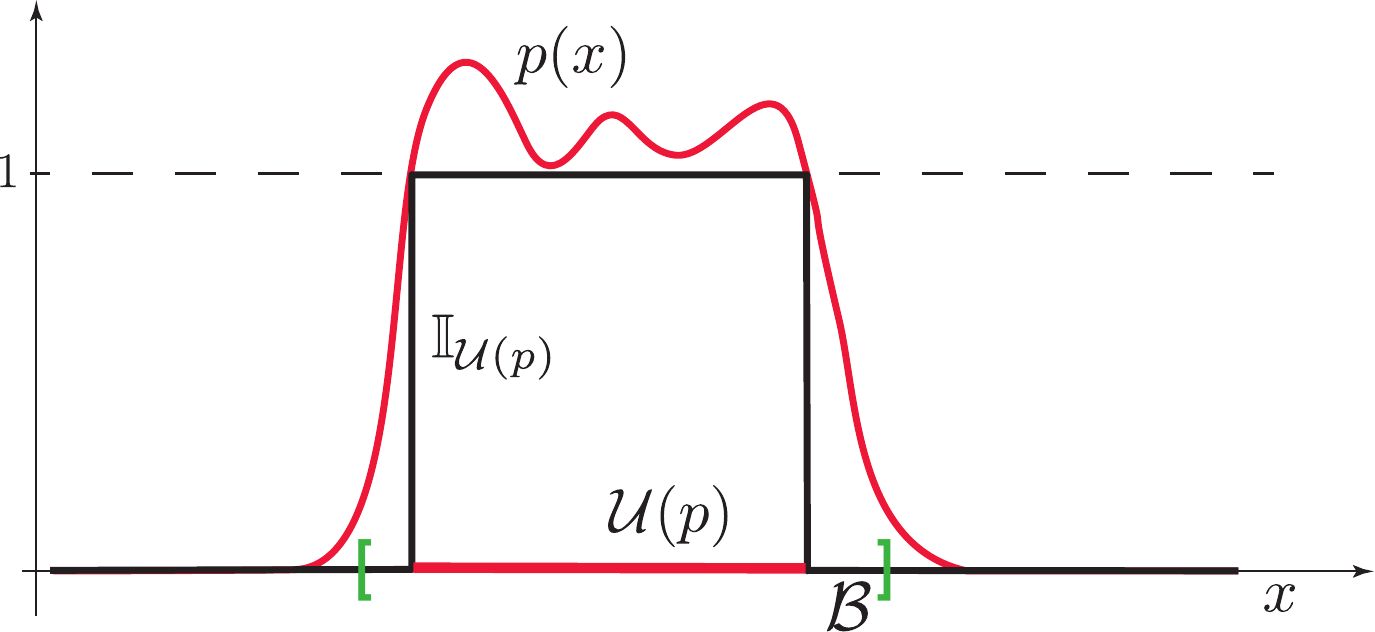}
}
\caption{Illustration of Chebychev's inequality: the polynomial is always greater or equal than the indicator function of $p(x)\ge 1$, hence the integral of $p$ over $\B$ is always an upper bound of the volume of $\U(p)$. \label{fig:Cheb}}
\end{figure}

{Assume now, without loss of generality, that the polynomial $p$ used to build the PSS is non-negative on $\B$.}
Then, observe that by definition of PSS  (see Figure \ref{fig:Cheb} for an illustration) we have
\[
p \geq \one_{\U(p)} \:\:\mathrm{on}\:\: \B.
\]

Hence, integrating both sides we get the following inequality
\begin{equation}\label{cheb0}
\int_\B p(x)dx \geq \int_\B \one_{\U(p)}(x)dx = \vol\:\U(p).
\end{equation}
This inequality is indeed  
widely used in probability, where it goes under the name of Chebyshev's inequality,
see e.g. \cite[\S 2.4.9]{AshDol:00}.
Note that, since the polynomial $p$ is nonnegative on $\B$, then 
the left-hand side of inequality \eqref{cheb0} corresponds to the $L^1$-norm of $p$ on $\B$, so that  the inequality simply becomes
\begin{equation}\label{cheb}
\|p\|_1\geq \vol\:\U(p).
\end{equation}

These derivations motivate us to the formulation of the following $L^1$-norm minimization problem,
which we choose as a surrogate of the original minimum volume outer PSS introduced in Problem 1.

\begin{problem}[Minimum {$L^1$}-norm outer {PSS}]
\label{minL1PSS}
Given a semialgebraic set $\K$, a bounding set $\B\supseteq \K$, and a degree $d$, 
 solve the optimization problem
\begin{equation}\label{l1}
\begin{array}{rcll}
w^*_d & \doteq & \displaystyle \inf_{p \in \Poly{d}} & \|p\|_1 \\
&& \mathrm{s.t.} & p \geq 0 \:\:\mathrm{on}\:\: \B \\
&&& p \geq 1 \:\:\mathrm{on}\:\: \K.
\end{array}
\end{equation}
\end{problem}

Note that a $L^1$-norm minimization
approach was originally proposed in \cite{HeLaSa:09}
for the numerical computation of the volume and of the higher order moments
of a semialgebraic set. The intuition underlying the formulation of Problem~\ref{minL1PSS}
is similar.
We now elaborate on some of the characteristics of the minimum $L^1$-norm outer PSS problem defined above.
First note that, for fixed $d$, when solving Problem~\ref{minL1PSS} we are minimizing an upper-bound on the volume of the PSS. Thus, the solution is expected to be a good approximation of the set $\K$.
Second, it can be shown that, as the degree $d$ increases, the Chebyshev bound~\eqref{cheb}
becomes increasingly tight. Indeed, the following fundamental result shows that the proposed solution converges to the minimum volume outer PSS.

\begin{theorem}\label{cvg}
{Given $d \in {\mathbb N}$, the infimum in problem (\ref{l1}) is attained for a polynomial
$p^*_d \in \Poly{d}$. Moreover, $w^*_d \geq v^*_d$ and $\U(p^*_d) \supseteq \K$.
Finally $w^*_d \geq w^*_{d+1}$ and $\lim_{d\to\infty} w^*_d = \lim_{d\to\infty} v^*_d = \vol\:\K$.}
\end{theorem}

\noindent
\textbf{Proof:}
{Let us first extend optimization problem (\ref{l1}) to continuous functions:
\begin{equation}\label{plp}
\begin{array}{rcll}
w^* & \doteq & \displaystyle \inf_f & \displaystyle \int_{\B} f(x)dx \\
&& \mathrm{s.t.} & f \in \mathcal{C}_+(\B) \\
&&& f-1 \in \mathcal{C}_+(\K)
\end{array}
\end{equation}
where $\mathcal{C}_+(\B)$ denotes the convex cone of non-negative continuous functions on $\B$. Observe
that since $f$ is non-negative on $\B$, the objective function $\|f\|_1 = \int_{\B} f(x)dx$ is linear.
Problem (\ref{plp}) is an infinite-dimensional linear programming (LP) problem in cones of non-negative continuous
functions. It has a dual LP, in infinite-dimensional dual cones of measures:
\begin{equation}\label{dlp}
\begin{array}{rcll}
v^* & \doteq & \displaystyle \sup_{\mu,\hat{\mu}} & \displaystyle \int \mu(dx) \\
&& \mathrm{s.t.} & \mu(dx)+\hat{\mu}(dx) = \one_{\B}(x)dx \\
&&& \hat{\mu} \in \mathcal{C}'_+(\B)\\
&&& \mu \in \mathcal{C}'_+(\K)\\
\end{array}
\end{equation}
where $\mathcal{C}'_+(\B)$ is the cone of non-negative continuous linear functionals on $\mathcal{C}_+(\B)$,
identified with the cone of Borel regular non-negative measures on $\B$, according to a Riesz
Representation Theorem \cite[Section 21.5]{RoyFit:10}. In LP (\ref{dlp}) the right hand side in the equation
is the Lebesgue measure on $\B$. Since the mass of non-negative measures $\mu$ and $\hat{\mu}$ is
bounded, it follows from Alaoglu's Theorem on weak-star compactness \cite[Section 15.1]{RoyFit:10} that
the supremum is attained in dual LP (\ref{dlp}) and that
there is no duality gap between the primal and dual LP, i.e. $v^*=w^*$, see also e.g. \cite[Theorem IV.7.2]{Barvinok:02}.

Moreover, as in the proof of \cite[Theorem 3.1]{HeLaSa:09}, it holds $v^*=\vol\:\K$. To see this, notice first
that the constraint $\mu+\hat{\mu}  = \one_{\B}$ jointly with $\mu \in \mathcal{C}'_+(\K)$ imply
that $\mu \leq  \one_{\K}$ and hence $\int \mu \leq \int \one_{\K} = \vol\:\K$ for every $\mu$ feasible in LP (\ref{dlp}).
In particular, this is true for an optimal $\mu^*$ attaining the supremum, showing $\int \mu^* = v^* \leq \vol\:\K$.
Conversely, the choice $\mu=\one_{\K}$ is trivially feasible for LP (\ref{dlp}) and hence suboptimal, showing
$v^* \geq \int \mu = \vol\:\K$. From this proof it also follows that the only optimal
solution to LP (\ref{dlp}) is the pair $(\mu^*,\hat{\mu}^*) = (\one_{\K},\one_{\B \backslash \K})$.

Now let us prove the statements of the Theorem:
\begin{itemize}
\item Attainment of the infimum in problem (\ref{l1}) follows from continuity (actually linearity) of the objective function
$\|p\|_1 = \int_{\B} p(x)dx = 0$ which is a norm (i.e. $\|p\|_1 = 0$ implies $p=0$ for $p \in \Poly{d}$)
and compactness of the set $\{p \in \Poly{d} \: :\: p \in \mathcal{C}(\B), \: \|p\|_1 \leq r\}$ for any fixed $r>0$.
\item $w^*_d \geq v^*_d$ follows readily from (\ref{cheb}).
\item $w^*_d \geq w^*_{d+1}$ follows readily from $\Poly{d} \subset \Poly{d+1}$.
\item Finally, $\lim_{d\to\infty} w^*_d = \vol\:\K$ is a consequence of $v^*=w^*=\vol\:\K$ (proven above)
and the Stone-Weierstrass Theorem  \cite[Section 12.3]{RoyFit:10} allowing to approximate uniformly on $\B$
by polynomials any continuous function in a minimizing sequence for LP (\ref{plp}), i. e. $\lim_{d\to\infty} w^*_d = w^*$.
\end{itemize}
}
$\Box$.

Some remarks are at hand regarding the above result, which represents one of the main contributions of the paper.

\begin{remark}[Convergence almost everywhere] 
\label{sec:trace}Note that Theorem \ref{cvg} implies that, for high enough order of approximation,  the PSS obtained by minimizing the $L^1$-norm of the polynomial defining it can be ``arbitrarily close'' to the semialgebraic set of interest. More precisely, as $d\rightarrow \infty$, $\|p^*_d\|_1$ and, as a consequence $\vol\:\U(p_d^*)$, converges to  $\vol\:\K$. Since,  $\K \subseteq \U(p_d^*)$, the Lebesgue measure of the difference between these sets converges to zero. In other words, one has almost everywhere convergence. From Theorems 2.5.1
and 2.5.3 in \cite{AshDol:00} the convergence is also almost
uniform, up to extracting a subsequence.
\end{remark}

\begin{remark}[Trace minimization] 
\label{sec:trace}
We provide a geometric interpretation
that further justifies the approximation of the minimum-volume PSS with the minimum $L^1$-norm {PSS}.
To this end, we first note that the objective function in \eqref{l1} reads
\begin{equation} \label{gram}
\|p\|_1= \int_\B p(x)dx = \int_\B \pi^T_\delta(x)P\pi_\delta(x)dx
= \trace\left(P\int_\B \pi_\delta(x)\pi^T_\delta(x)dx\right)
= \trace\:PM
\end{equation}
where 
\[
M\doteq\int_\B \pi_\delta(x)\pi^T_\delta(x)dx
\]
is the matrix of moments of the Lebesgue measure on $\B$ in the basis $\pi_\delta(x)$.
Note that,  if the basis in equation \eqref{gram} is chosen such that its entries
are orthonormal with respect to the (scalar product induced by the)
Lebesgue measure on $\B$, then $M$ is the
identity matrix and inequality (\ref{cheb}) becomes
\[
\trace\:P \geq \vol\:\U(p)
\]
which indicates that, under the above constraints,
minimizing the trace of the Gram matrix $P$
entails minimizing the volume of $\U(p)$.
It is important to remark that, in the case of quadratic polynomials, i.e.\ $d=2$, we retrieve
the classical trace heuristic used for volume minimization of ellipsoids, 
see e.g. \cite{DuPoWa:96}. 
Indeed, if $\B=[-1,1]^n$, then the basis $\pi_1(x)=\frac{\sqrt{6}}{2}x$
is orthonormal with respect to the Lebesgue measure on $\B$ and
$\|p\|_1 = \frac{3}{2}\trace\:P$. Moreover, note that the constraint that $p$ is
nonnegative on $\B$ implies that the curvature of the boundary
of $\U(p)$ is nonnegative, hence that $\U(p)$ is convex. Thus, $\U(p)$ is indeed
an ellipsoid. 
\end{remark}
\vskip .1in

\begin{remark}[Choice of $\B$] 
\label{sec:trace}
We finally remark that, as previously noted, the assumption of the bounding set $\B$ being an hyperrectangle can be easily relaxed. Indeed, in order to develop a computationally manageable optimization in Problem 2,  $\B$ can be selected as a semialgebraic set, provided that the polynomials defining the set should be such that the objective function in problem
(\ref{l1}) is easy to compute. In particular, if
\[
p(x)=\pi^T_d(x) p = \sum_\alpha p_\alpha [\pi_d(x)]_\alpha
\]
then
\[
\int_\B p(x)dx = \sum_\alpha p_\alpha \int_\B [\pi_d(x)]_{\alpha}dx
= \sum_\alpha p_\alpha y_\alpha
\]
and we should be able to compute easily the moments
$
 \int_\B [\pi_d(x)]_{\alpha} dx
$
of the Lebesgue measure on $\B$ with respect to the basis $\pi_d(x)$.
\end{remark}

\section{{LMI hierarchy to compute the PSS}} 
\label{sec:sos}
 
In this section, we provide the basic details on the numerical  computation of the solution of the minimum $L^{1}$-norm PSS introduced in Problem 2.
Note that, in problem~\eqref{l1}, we aim at finding a polynomial $p \in \Poly{d}$ such that i) $p$ is positive on $\B$, and ii) $p-1$ is positive on $\K$. In order to obtain a numerically solvable problem, we {enforce} positivity  by requiring the polynomial to be SOS, and use Putinar's Positivstellensatz; e.g., see  \cite{Putinar:93,Lasserre:01,CGTV:03,Parrilo:03}. More precisely,  fix $r \in {\mathbb N}$, and consider the problem
 \begin{align}
{w^*_{2r,d}}  = &  \min_{p\in \Poly{d}} \int_{\B} p(x) dx   \label{l1_sos} \\
              &\text{s.t.} \nonumber\\ 
             & \left.\begin{aligned}
              & p(x) = s_{0,\B}(x) + \sum_{j=1}^{n} s_{j,\B}(x) (x_j-a_j)(b_j-x_j) \notag \\
                            &       s_{0,\B} \in \SoS{2r} \notag \\
                            & s_{j,\B} \in \SoS{2(r-1)}, \quad j=1,2,\ldots,n \notag\\
       \end{aligned}
 \right\}\ \ p(x) \text{ positive on } \B {=[a,b]}\\
             & \left.\begin{aligned}
              & p(x)-1 = s_{0,\K}(x) + \sum_{i=1}^{m} s_{i,\K}(x) g_i(x) \qquad\qquad\notag \\
              &       s_{0,\K} \in \SoS{2r} \notag \\
              & s_{i,\K} \in \SoS{2(r-r_i)}, \quad i=1,2,\ldots,m. \notag 
       \end{aligned}
 \right\}\ \ p(x)-1 \text{ positive on } \K \\ \nonumber
  \end{align}
{where $r_i$ is the smallest integer greater than half the degree of $g_i$ for $i=1,2,\ldots,m$.}  
It should be noted  that the objective function of problem \eqref{l1_sos} is an easily computable linear function of the coefficients of the polynomial $p$. Moreover,  the constraints can be recast in terms of Linear Matrix Inequalities (LMIs); see, for instance, {\cite{Lasserre:01}}). Several \textsc{Matlab} toolboxes have efficient and easy to use {interfaces to model} problems of the form above; e.g., see YALMIP~\cite{yalmip}. 

{Not only we can numerically solve problem~\eqref{l1_sos}, but the following result holds. This theorem  is an immediate consequence of the results in~\cite{Putinar:93}.

\begin{theorem}
Let us denote by $p^*_{2r,d}$ a solution of problem ~\eqref{l1_sos}. Then, the following hold
 \begin{enumerate}[i)]
 	\item for each $d \in {\mathbb N}$, the value of problem~\eqref{l1_sos} converges to the value of problem~\eqref{l1} as $r\to\infty$, {i.e. $\lim_{r\to\infty} w^*_{2r,d} = w^*_d$},
 	
 	\item for any $2r\geq d$, $p^*_{2r,d} \geq 0$ {on $\B$},
 	\item for any $2r\geq d$, $p^*_{2r,d} \geq 1$ {on $\K$}.
 \end{enumerate}
\end{theorem}}

We conclude that $p_{2r,d}^*$ can be used to compute a PSS approximation for $\K$. 
For our numerical examples,
we have used the YALMIP~\cite{yalmip} interface for Matlab  to model the
LMI optimization problem (\ref{l1_sos}) and the SDP solver SeDuMi~\cite{SEDUMI}
to numerically solve the problem. Since the degrees of
the semialgebraic sets we compute are typically low (say less than 20),
we did not attempt to use {alternative} polynomial bases (e.g. Chebyshev
polynomials) to improve the quality and
resolution of the optimization problems; see  \cite{HeLaSa:09}
for a discussion on these numerical matters {in the context of semialgebraic set volume approximation.}

\subsection{Computing Bounding Box $\B$}
\label{sec:Bab}

As noted in \cite[Remark 1]{CePiRe:12}, an outer-bounding hyper-rectangle $\B=[a,b]$ of a given semialgebraic set $\K$ can be found
by solving relaxations of the following polynomial optimization problems
\begin{eqnarray*}
	a_{j} = \arg\min_{x\in\mathbb{R}^{n}} x_{j} \text{ subject to } x \in\K,\quad j = 1,...,n,\\
	b_{j} = \arg\max_{x\in\mathbb{R}^{n}} x_{j} \text{ subject to } x \in\K,\quad j = 1,...,n,
\end{eqnarray*}
which compute the minimum and maximum value of each component of the vector $x$ over the semialgebraic set $\K$.

To illustrate how this can be done, let us concentrate on approximating the value of $a_j$. First, note that the problem of computing $a_j$ is equivalent to solving the following polynomial optimization problem
\[
a_j= \max y\, \text{ subject to } x_j - y\geq 0 \text{ for all } x \in \K.
\]
Then, formulate the following convex optimization problem
\begin{align*}
a_{j,2r}  = &  \max y \\
s. t. \ \ \ &  x_j-y= s_0(x) + \sum_{i=1}^{m} s_i(x) g_i(x) \\
&       s_0 \in \SoS{2r}; \\
& s_i \in \SoS{2(r-r_i)}; \quad i=1,2,\ldots,m.
\end{align*}
Using the same reasoning as above, it can be shown that: i) $a_{j,2r} \leq a_j $  for all $r$, and ii) {$\lim_{r\to\infty} a_{j,2r} = a_j$}. Moreover, the problem above can be recast  as an {LMI} optimization problem. 

\section{Inner approximations} \label{sec:inner}

The approach described in the previous sections can be readily extended to derive \textit{inner} approximations of the set $\K$, in the spirit of \cite{DaGaPo:02,HeLa:12,HenLou:12}.
The idea is just to construct an optimal outer PSS of the \textit{complement set}
\begin{align*}
\overline \K \doteq \B \setminus \K & =
\{x \in {\mathbb R}^n : g_1(x) < 0\text{ or }\cdots\text{ or } g_m(x) \le 0, \:i=1,2,\ldots,m\} \cap \B\\ \nonumber \\
& = \left( \K_{1}\cup\K_{2}\cup\cdots\cup\K_{m}\right) \cap \B, \nonumber
\end{align*}
with $\K_{j}\doteq
\{x \in {\mathbb R}^n : g_j(x) < 0\}$.


Note that, since the set whose indicator function we want to approximate is a union of basic semialgebraic sets, the $L^{1}$ optimization problem to be solved becomes
\begin{equation}\label{l1-inner}
\begin{array}{ll}
 \displaystyle \min_{p \in \Poly{d}} & \|p\|_1 \\
 \mathrm{s.t.} & p \geq 0 \:\:\mathrm{on}\:\: \B \\
& p \geq 1 \:\:\mathrm{on}\:\: \K_{1}\\
& p \geq 1 \:\:\mathrm{on}\:\: \K_{2}\\
& \vdots \\
& p \geq 1 \:\:\mathrm{on}\:\: \K_{m}\\
\end{array}
\end{equation}
{and let $p^*_d$ attain the minimum}.
The corresponding optimal inner approximation is given by the polynomial sublevel set
\[
\V(p^*_d) \doteq \{ x \in \B : p^*_d(x) \leq 1 \}.
\]
In this case, one can think of the polynomial $1-p^*_d$ as a lower bound for the indicator function of the set $\K$.  

Given the fact that the optimization problem~\eqref{l1-inner} provides an outer approximation of the set~$\overline{\K}$, one has the following result whose proof is similar to that of Theorem~\ref{cvg}.

\begin{corollary}
{For all $d \in {\mathbb N}$} it holds $\V(p_d^*) \subseteq \K$. {Moreover $\lim_{d\to\infty} \vol\:V(p_d^*) = \vol\:\K$}.
\end{corollary}

As before, to be able to numerically approximate the solution of problem~\eqref{l1-inner}, we replace its {polynomial
positivity} constraints by their  {LMI} approximations as  described in Section~\ref{sec:sos}.

\section{Numerical examples} \label{sec:examples1}

In this section, we present several examples that illustrate the performance of the proposed approach.
\subsection{{Discrete-time} stabilizability region}
\label{example-didier}
\begin{figure}[!ht]
\centerline{
\includegraphics[width=9cm]{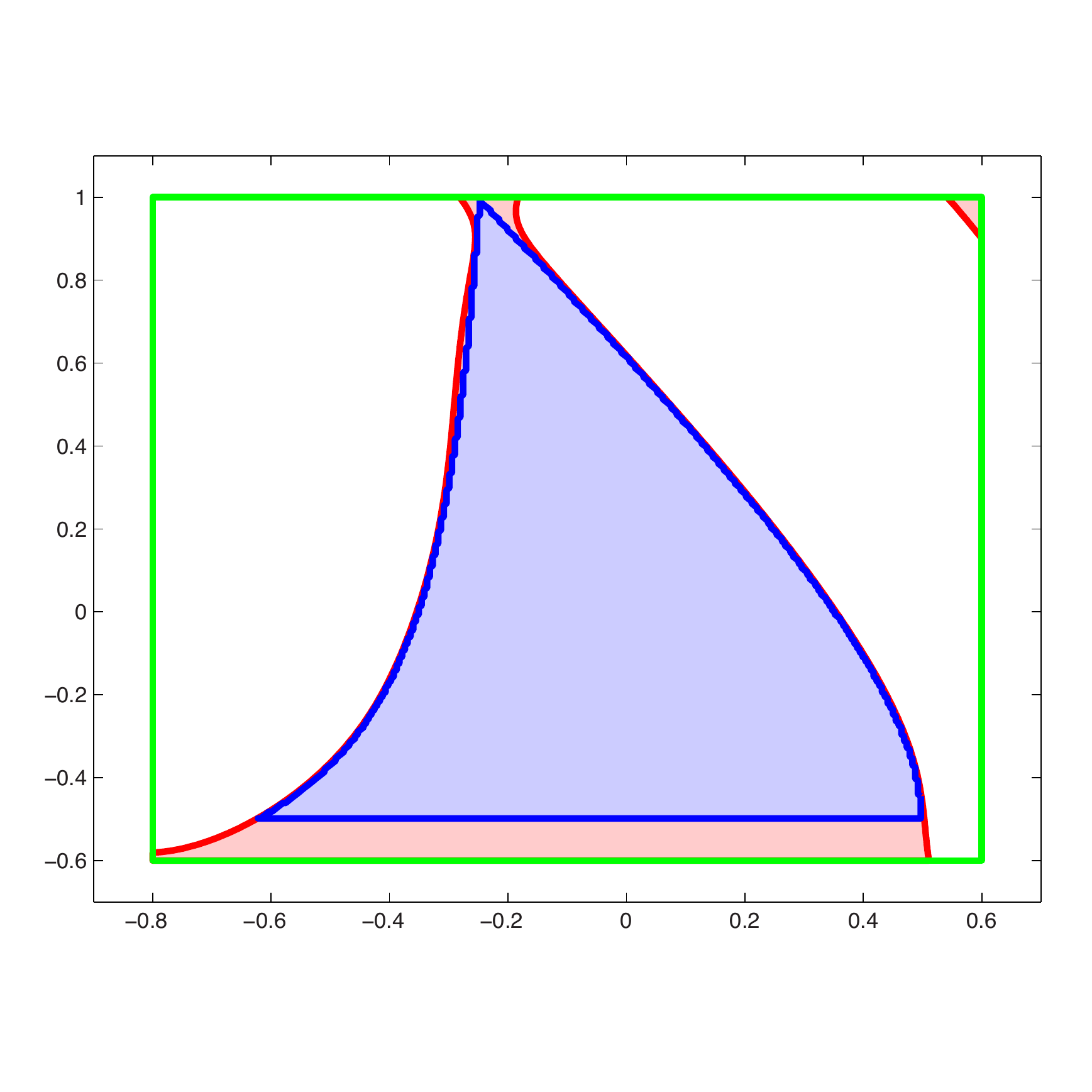}
\includegraphics[width=9cm]{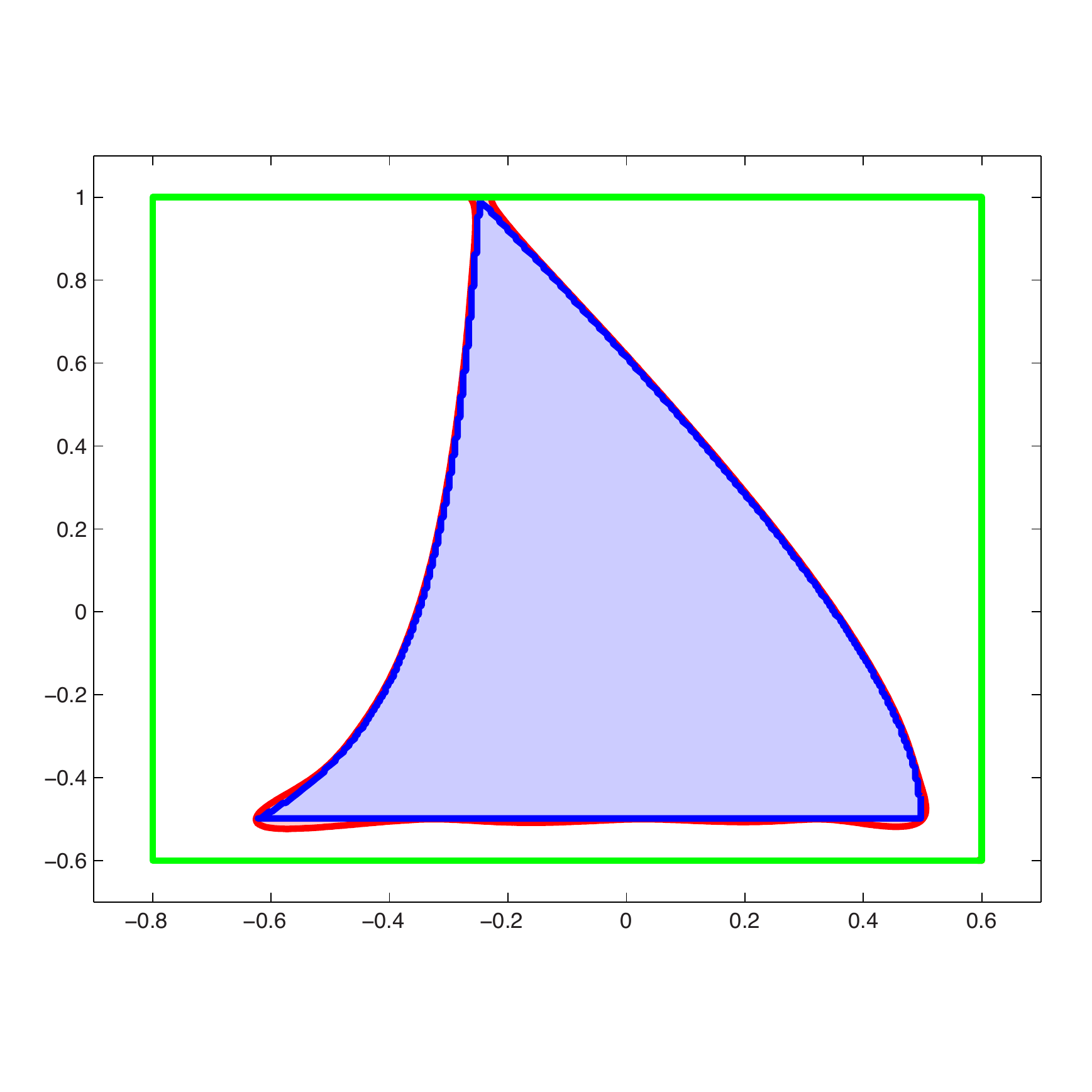}
}
\caption{Degree 6 and degree 20 outer PSS approximation (red)  
of {stabilizability} region {$\K$} (inner {surface} in light blue). The green box corresponds to the bounding set $\B$. \label{fig:schur4outer}}
\end{figure}

\begin{figure}[!ht]
\centerline{
\includegraphics[width=11cm]{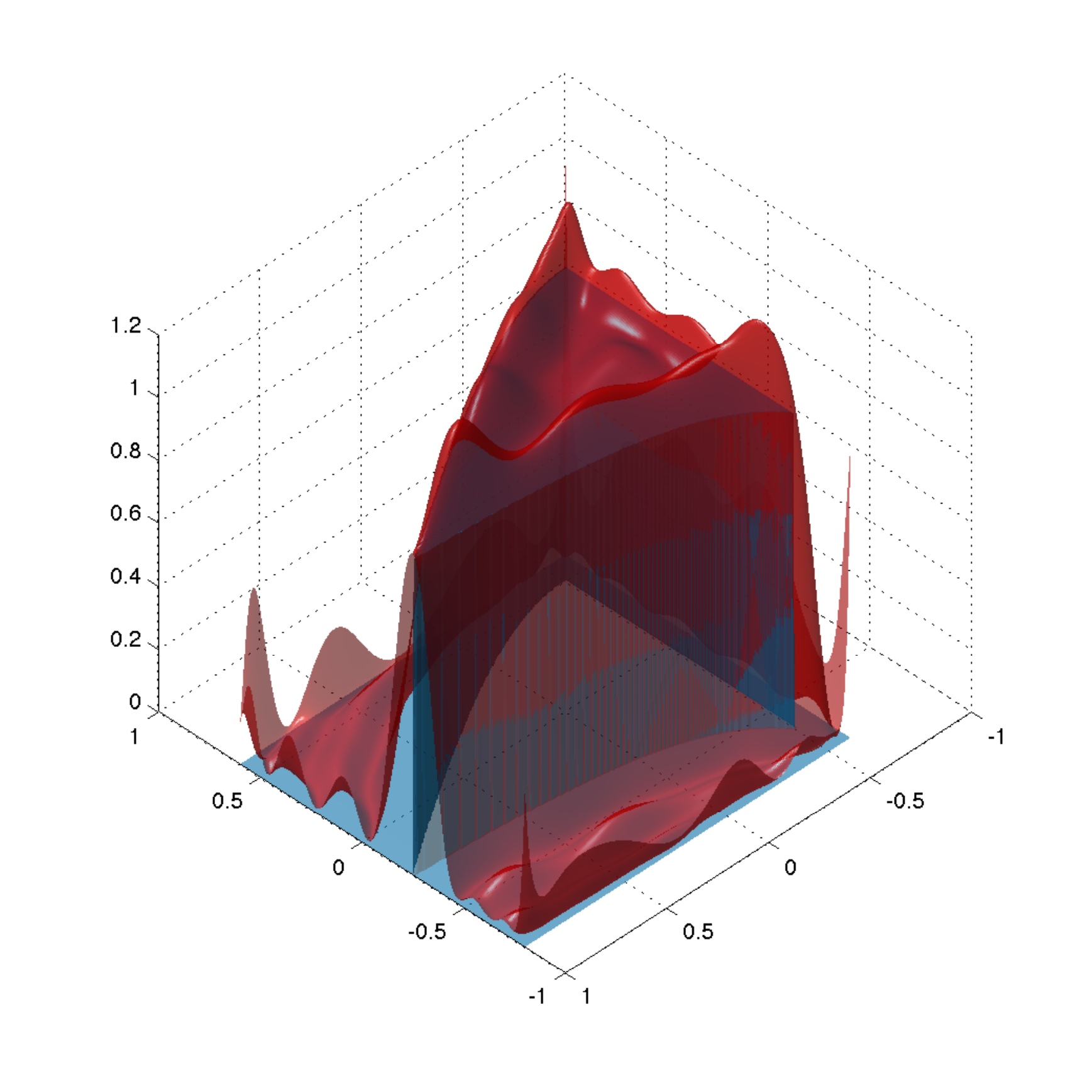}
}
\caption{Degree 20 polynomial approximation ({upper} surface in red)
of the indicator function ({lower} surface in blue) of the nonconvex planar stabilizability
region {$\K$}.\label{fig:schur4}}
\end{figure}

As a  control-oriented illustration of the {PSS} approximation described in this paper, consider
\cite[Example 4.4]{HeLa:12} which is a degree 4 discrete-time polynomial $z \in {\mathbb C} \mapsto 
x_2+2x_1z-(2x_1+x_2)z^3+z^4$
to be stabilized by means of 2 real control parameters $x_1, x_2$.
In other words, we are interested in approximating the set $\K$
of values of $x_1, x_2$ such that this polynomial has its roots with
modulus less than one. An explicit basic semialgebraic description
of the stabilizability region is built using {the Schur stability criterion}, resulting in the following basic semialgebraic set:
\begin{eqnarray}
\label{K-didier}
\K &= \{x \in {\mathbb R}^2 \: :\: &g_1(x)=1+2x_2\geq 0, \\
&&g_2(x) =2-4x_1-3x_2 \geq 0,\nonumber \\
&& g_3(x) =10-28x_1-5x_2-24x_1x_2-18x^2_2 \geq 0,\nonumber \\
&& g_4(x) = 1-x_2-8x_1^2-2x_1x_2-x_2^2-8x_1^2x_2-6x_1x_2^2 \geq 0\}. \nonumber 
\end{eqnarray}
This set is nonconvex and it is included in the box
$\B=[-0.8,0.6]\times[-0.5,1.0]$. In Figure \ref{fig:schur4outer}  we represent the 
PSS {outer} approximations {of $\K$} for $d=6$ and $d=12$ respectively, while Figure \ref{fig:schur4}
shows the graph of the degree $d=12$ polynomial $p^*_{12,12}(x)$ constructed by solving optimization problem \eqref{l1_sos} with $2r=d$. 

\begin{figure}[!ht]
\centerline{
\includegraphics[width=9cm]{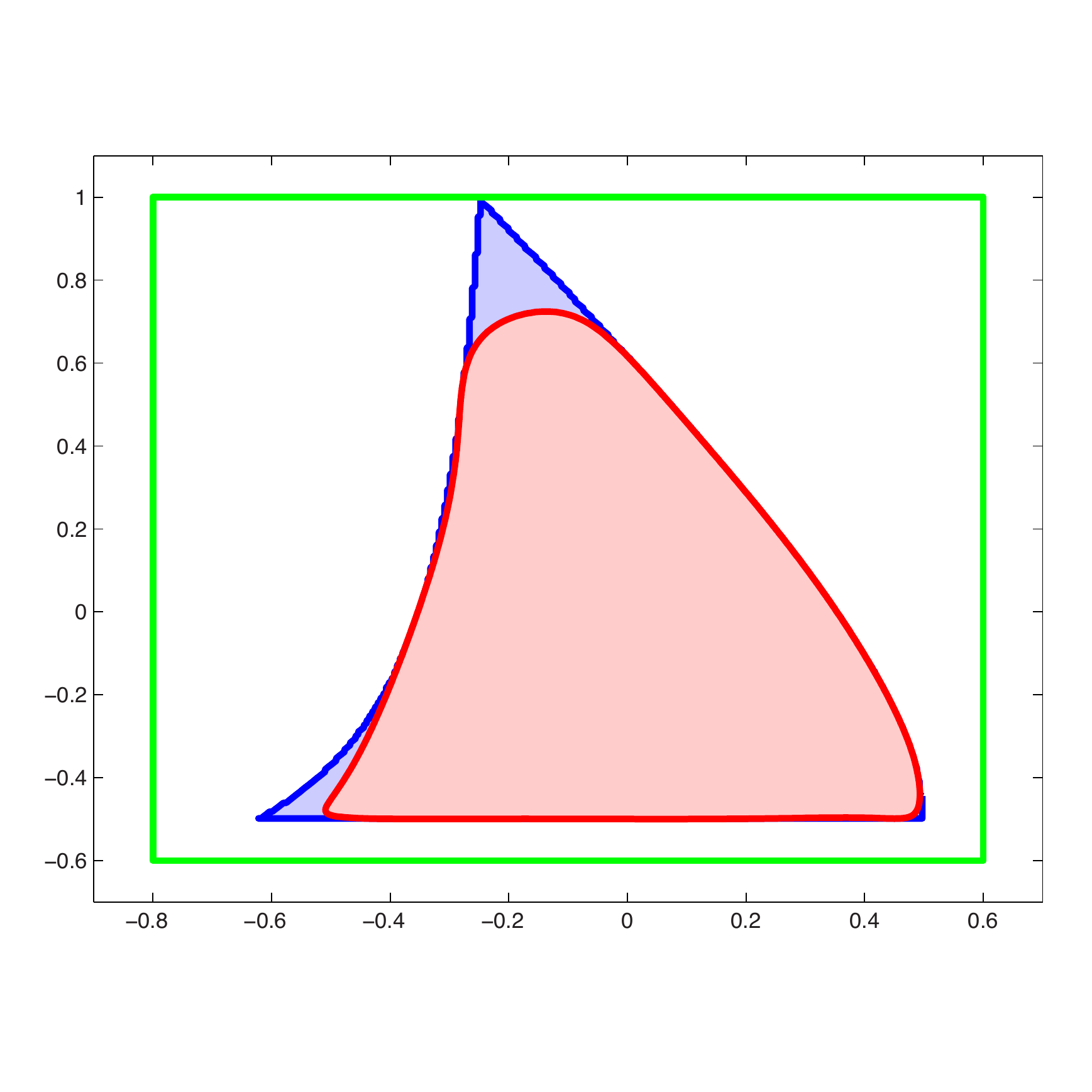}
\includegraphics[width=9cm]{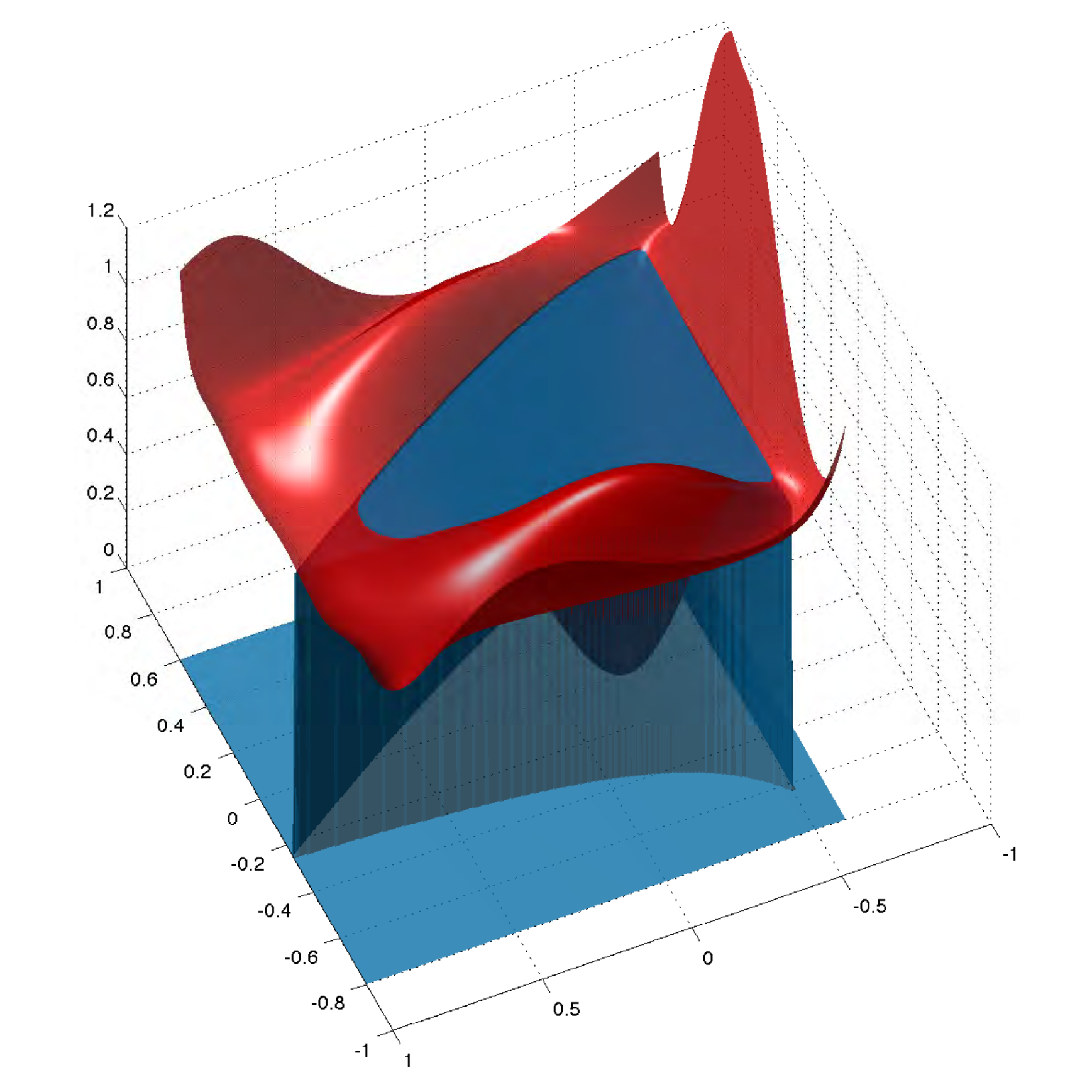}
}
\caption{{Left:} degree 8 inner PSS approximation (red)  
of {stabilizability} region {$\K$} (inner {surface} in light blue). {Right: degree 8 polynomial approximation (upper surface in red) of the indicator function (lower surface in blue} of $\K$. \label{fig:schur4inner}}
\end{figure}

As discussed before, we can also use the approach proposed in this paper to obtain inner approximations of {$\K$.}  In Figure~\ref{fig:schur4inner}, we depict the inner approximation obtained using optimization problem~\eqref{l1-inner} with $2r=d=8$. 

\subsection{{PID} stabilizability region}
\label{example-bhatt}
We now turn our attention to an example related to  fixed order controller design. Consider  \cite[Example 2.2]{BhDaKe:09}, in which the authors examine the problem of stabilizing the plant $P(s)=\frac{N(s)}{D(s)}$ where
\begin{eqnarray*}
N(s)&=&s^{3}-2s^{2}-s-1;\\
D(s)&=&s^{6}+2s^{5}+32s^{4}+26s^{3}+65s^{2}-8s+1.
\end{eqnarray*}
by means of a PID controller of the form $K\ped{PID}(s)=k\ped{P}+\frac{k\ped{I}}{s}+k_{D}s$.
In particular, they are interested in finding the set of stabilizing PID gains, that is the set of gains for which 
the closed-loop characteristic polynomial $sD(s)+(k\ped{I}+k\ped{P}s+k\ped{D}s^{2})N(s)$
is Hurwitz. For this special class of controllers, the authors provide a method based on the so-called signature of a set {of} properly constructed polynomials to determine the set of all PID gains that stabilize the plant. One should note that this procedure is not easily generalizable to more general classes of fixed order controllers.

\begin{figure}[!ht]
\centerline{
\includegraphics[width=9cm]{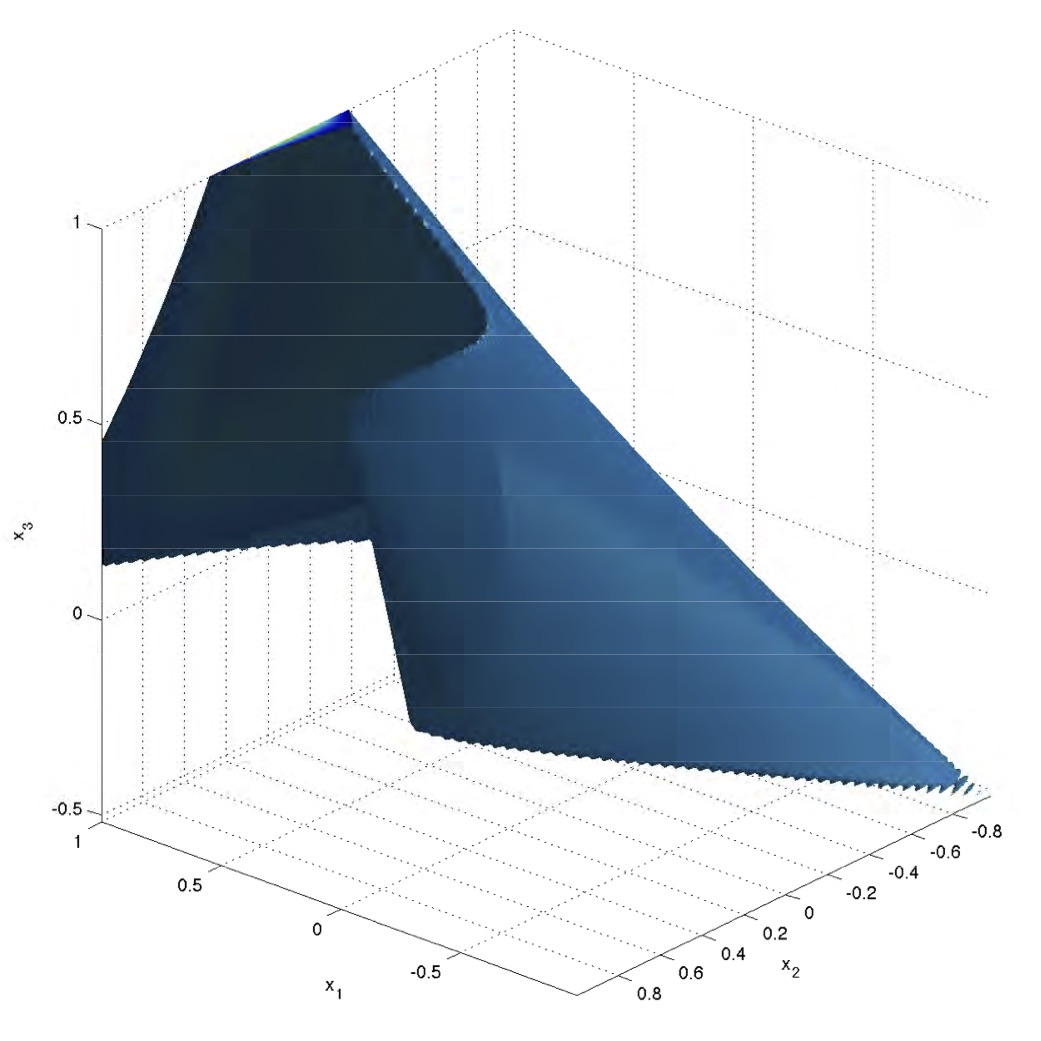}
\includegraphics[width=9cm]{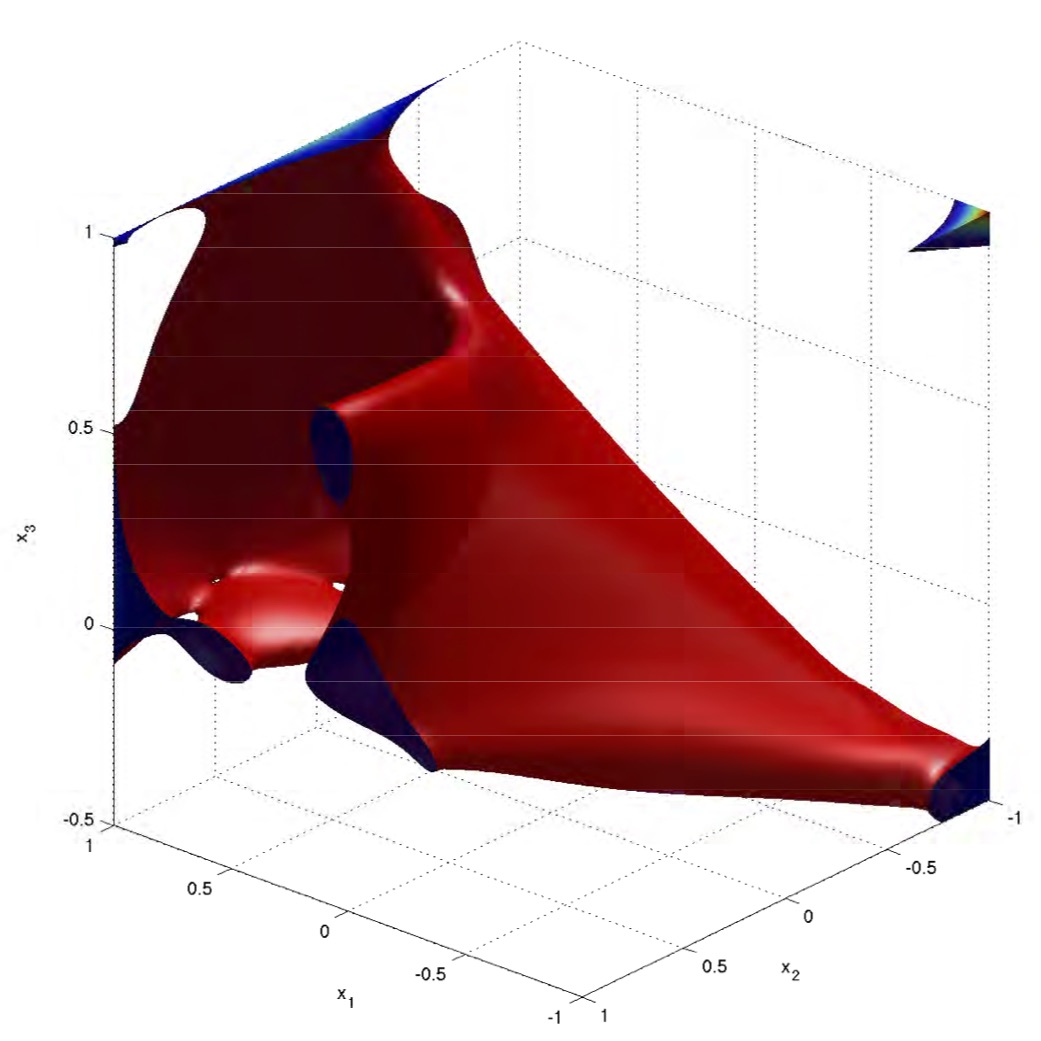}
}
\caption{{Left:} set of stabilizing {PID} gains.  {Right: its} degree 14 optimal outer PSS approximation.\label{fig:bhatt}}
\end{figure}

In our setup, we are interested in approximating the set 
\[
\K = \{x \in {\mathbb R}^3 \: :\: sD(s)+(k\ped{I}+k\ped{P}s+k\ped{D}s^{2})N(s)\text{ is Hurwitz},\;
k\ped{I}=25(x_{1}-1),\;
k\ped{P}=10(x_{2}-1.5),\;
k\ped{D}=10(x_{3}-1)\}
\]
with bounding box $\B=[-1,1]^{3}$. As one can see in Figure~\ref{fig:bhatt}, the approached proposed in this paper provides a very good approximation of the set of stabilizing gains, even for a PSS of relatively low order ($d=14$).

\section{Reconstructing/approximating Sets from a finite number of samples} \label{sec:finite_samples}

A particularly interesting case is when the semialgebraic set $\K$ is \textit{discrete}, that is, it consists of the union of $N$ points 
\[
\K = {\bigcup_{i=1}^N \{x^{(i)}\} \subset {\mathbb R}^n.}
\]
This situation arises for instance when the objective is to try to approximate a given set (possibly non-semialgebraic) from a given number of points in its interior. An example of this is the reconstruction of reachable sets by using randomly generated trajectories. This  setup
is discussed in \cite{DaLaSh:10,DHLS:15,HwStTo:03}. 

From a computation viewpoint, an important feature is that, in the case of a discrete set, the inclusion constraint $\K \subseteq \U(p)$ is equivalent to a finite number of inequalities 
\[
p(x^{(i)}) \geq 1, \quad i=1,\ldots,N
\]
which are \textit{linear} in the coefficients of $p$. This fact allows to deal with problems with rather large $N$. Moreover, in this latter case, where the number of points $N$ is large while the dimension $n$ is relatively small, the constraint that $p$ is nonnegative on $\B$ can also be (approximately) handled by linear inequalities 
\[
p(z^{(j)}) \geq 0, \quad j=1,\ldots,M
\]
enforced at a dense grid of points $z^{(j) }\in \B$, for $M$
sufficiently large. Hence, in this case one can construct a pure linear programming (LP) approach.
Note that, even if this approach does not guarantee that $p$ is nonnegative everywhere on $\B$, it still
ensures that $\K \subseteq \U(p)$, which is what matters
primarily in our approach.

\begin{figure*}[htb]
\centering
\includegraphics[width=58mm]{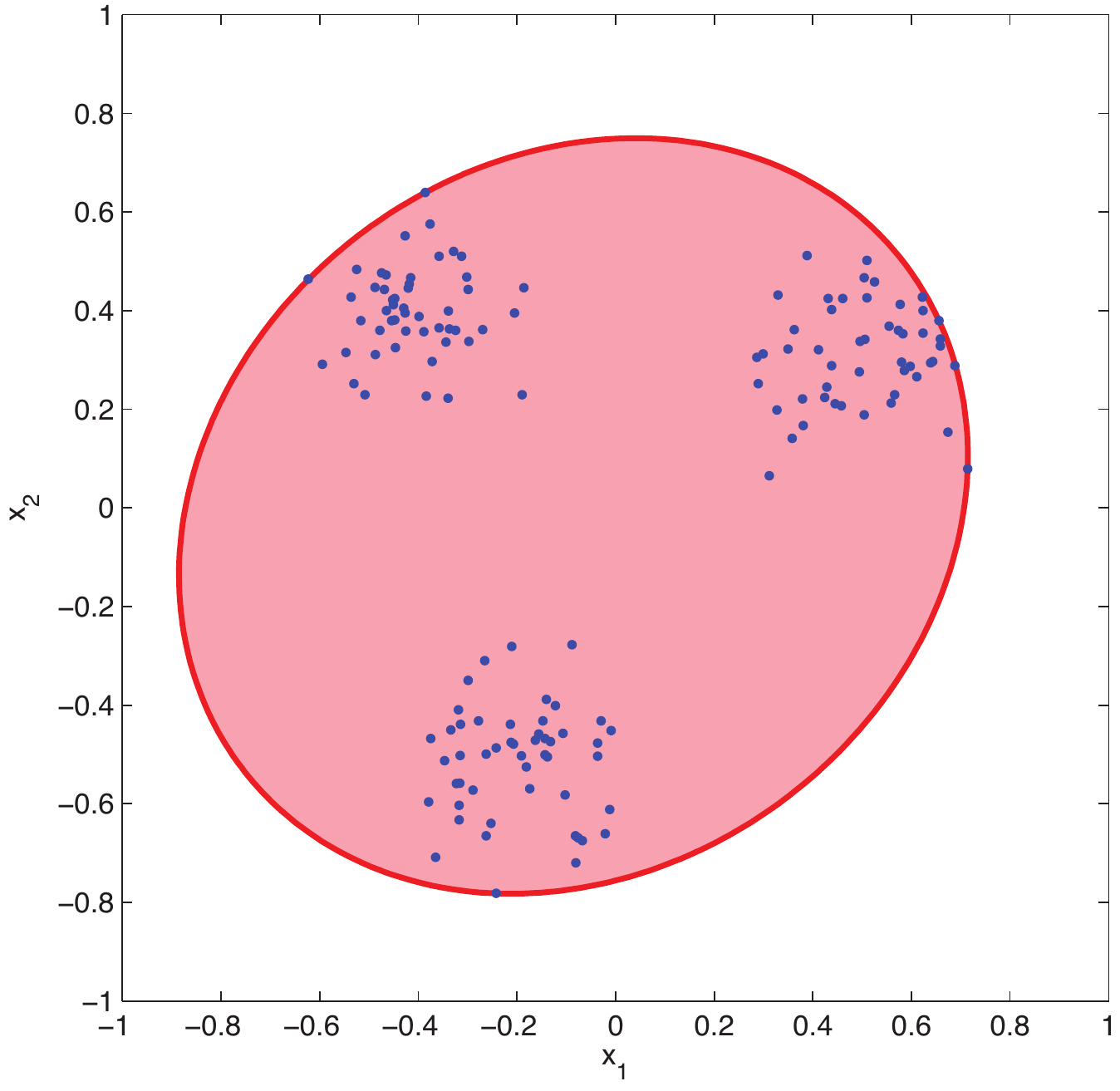}
\includegraphics[width=58mm]{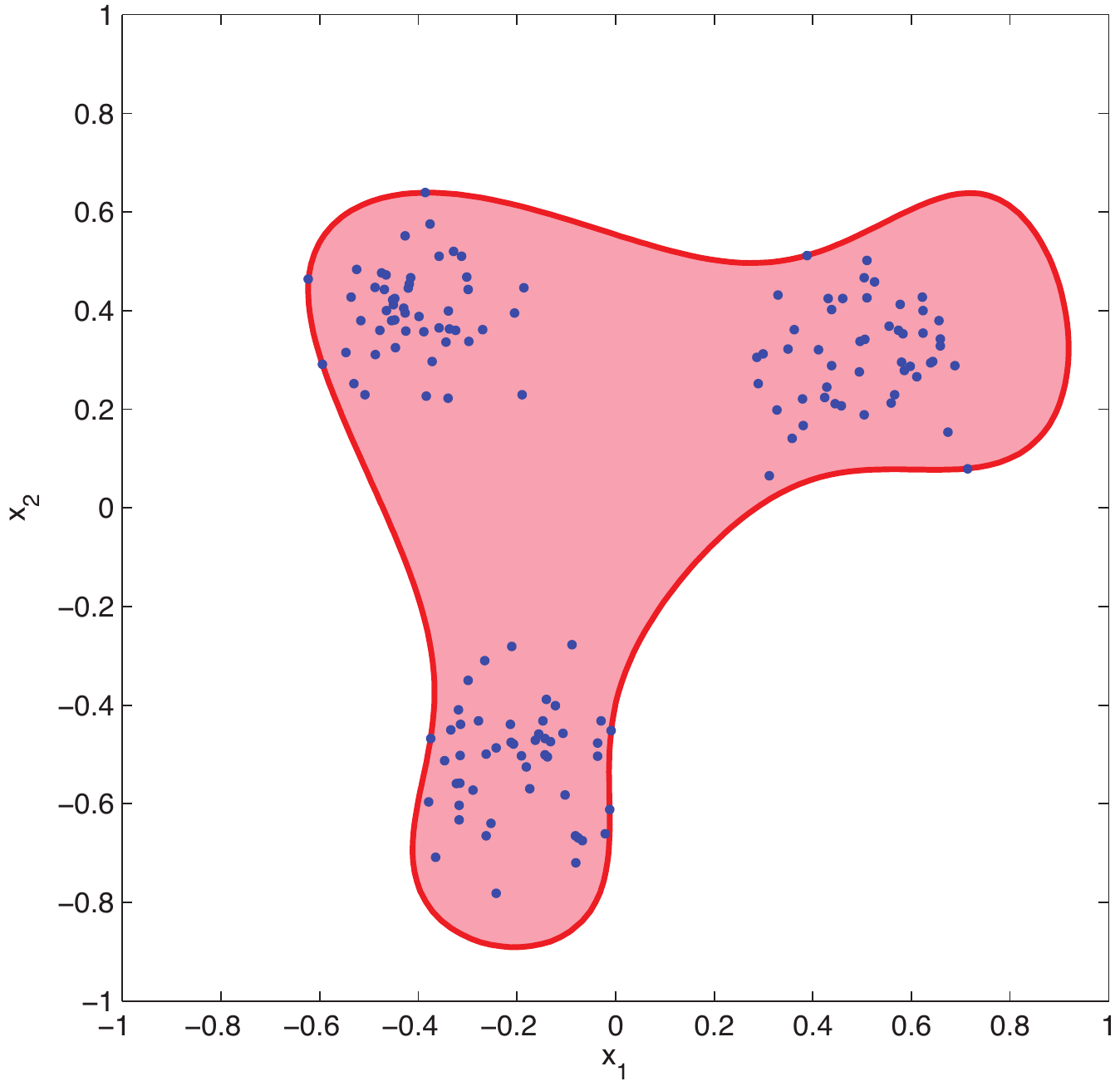}
\includegraphics[width=58mm]{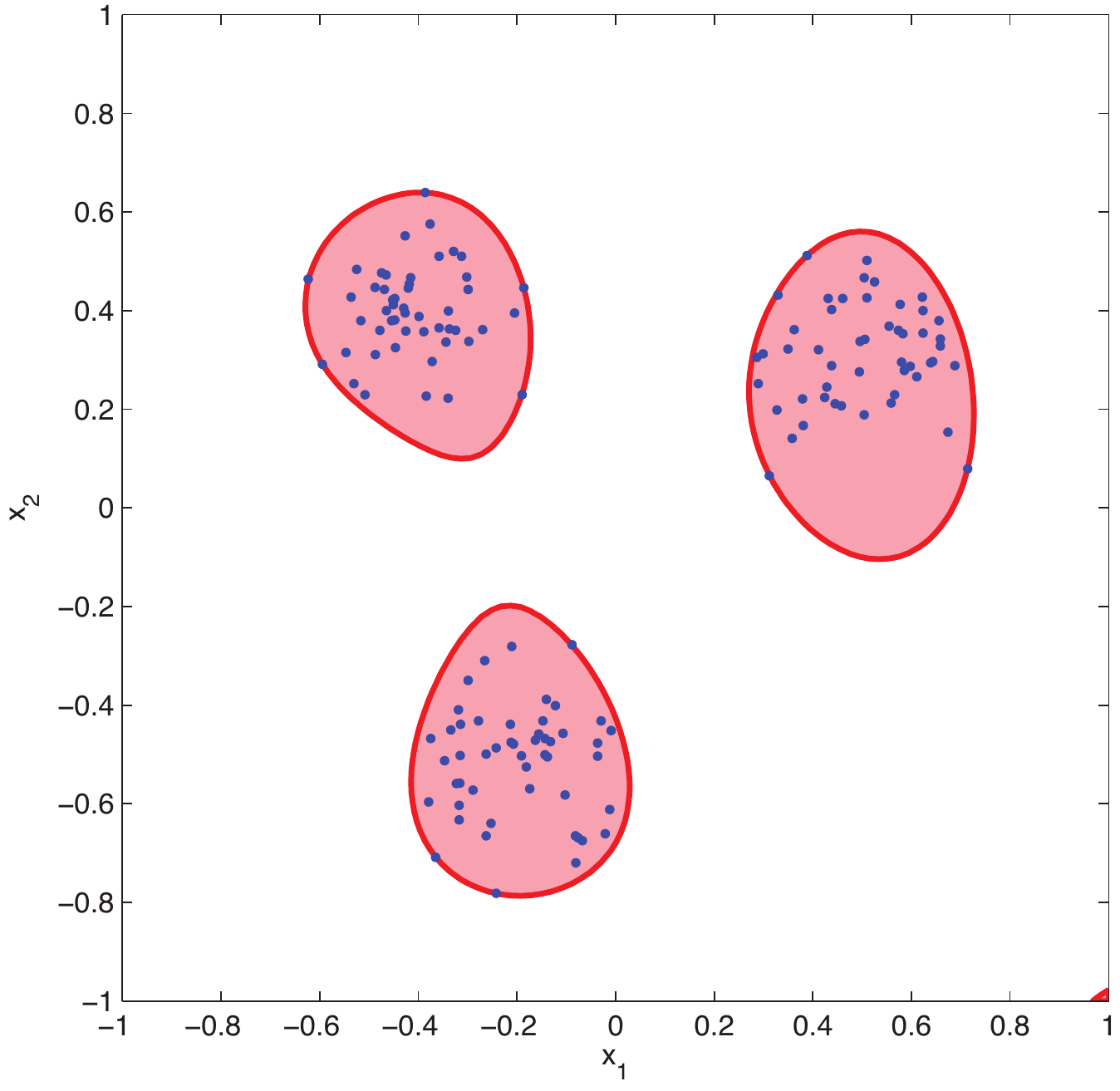}
\caption{Minimum $L^1$-norm PSS at 100 points (blue), for degree 2 (left),
5 (center), and 9 (right).\label{plane}}
\end{figure*}

To illustrate the performance of the proposed method, we first consider $N=100$ points in the box $\B=[-1,\:1]^2$. The points are generated mapping Gaussian points with variance $0.1I$ and mean value chosen with equal probability between $[0.4,\,0.3]^{T}$, $[-03,\,-0.5]^{T}$, $[-0.5,\,0.4]^{T}$. 
On Figure \ref{plane} we represent the solutions $p$ of degrees 2, 5, and 9 
of minimization problem (\ref{l1}). A few comments about the obtained solution are at hand. First, we see that the solution for $d=2$ 
corresponds to the L\"owner-John ellipsoid, see e.g. \cite[\S 4.9]{BenNem:01}.
Second, it can be observed that, as the degree of $p$ increases, 
the set $\U(p)$ becomes disconnected, so as to better capture the different regions where the points are concentrated. We note that, in the case of discrete points, it is not advisable to select high values of $d$, since  
indeed, in the limit, the optimal polynomial would correspond to a function with spikes corresponding to the location of the considered points. Finally, we remark that the {possible} side effects near the border of $\B$ on the right
{hand side} figure can be removed by enlarging the bounding set $\B$.

\begin{figure*}[htb]
\centering
\includegraphics[height=60mm]{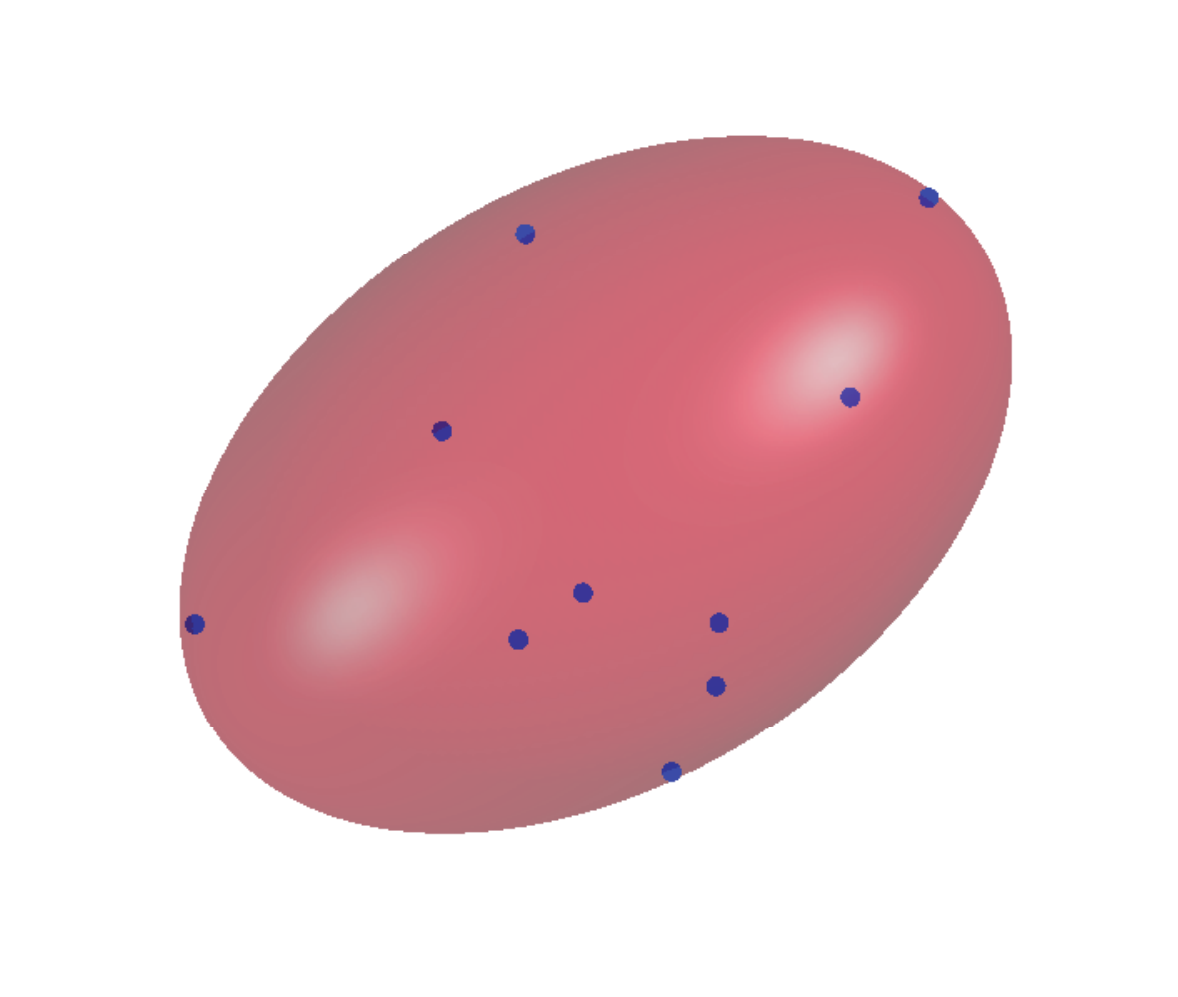}\qquad
\includegraphics[height=58mm]{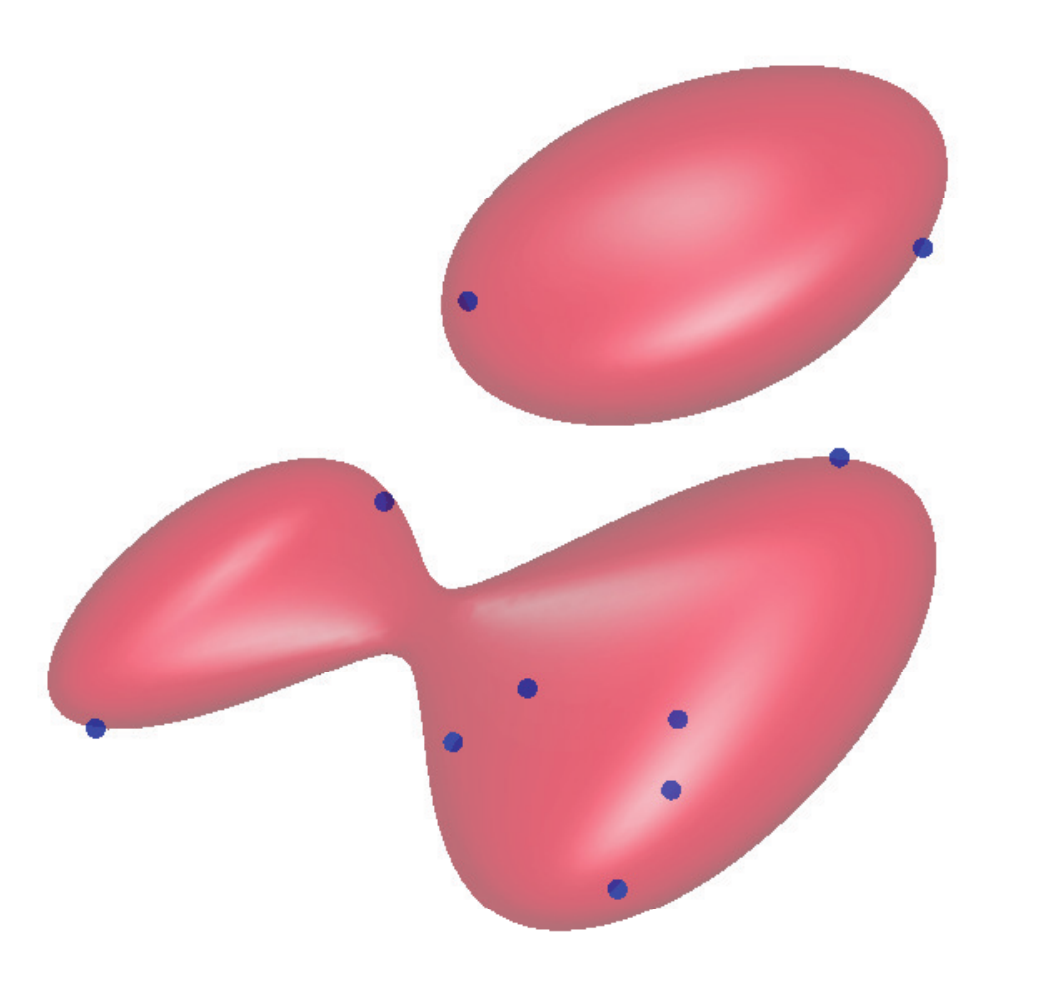}
\caption{Including the same 10 space points (blue) in PSS 
of degree 4 (left) and 10 (right).
\label{space}}
\end{figure*}

As a second illustrative example,  we consider $N=10$ points in $\B=[-1,\:1]^3$.
The solutions $p$ of degrees 4, 6, 9, and 14 of minimization problem (\ref{l1})
is depicted in Figure \ref{space}. Here too we observe that increasing the
degree of $p$ allows to capture point clusters in distinct connected components.

\section{Uniform sampling over semialgebraic sets} \label{sec:sampling}

In this section, we consider a problem that can be seen as the ``dual'' of the one considered in the previous section; that is, instead of trying to reconstruct/approximate {the indicator function of} an unknown set from points belonging to its interior, we  aim at developing systematic procedures for generating \textit{uniformly distributed samples} in a given semialgebraic set. 
This is an important problem since many system specifications lead to sets with a (complex) closed-form description, and being able to draw samples from these type of sets provides the means for the design of systems with a complex set of specifications.
In particular, the algorithm presented in this {section}
can be used to generate  uniform samples in the solution set of LMIs.

As before, we assume that the set of interest is a compact basic semialgebraic set defined as in \eqref{K-set},
and that there exists a bounding  hyper-rectangle $\B=[a,b]$ of the form \eqref{Bab}. Then, the problem we discuss in this section is the following.\\

\noindent
\begin{problem}[Uniform Sample Generation over $\K$] {\it Given a semialgebraic  set $\K$ defined in (\ref{K-set}) of nonzero volume,
generate $N$ independent identically distributed (i.i.d.) random samples $x^{(1)},\ldots,x^{(N)}$ uniformly distributed in~$\K$.}
\end{problem}

Let us start by describing the approach proposed to solve this problem. First, we define the uniform density over the set $\K$  as follows
\begin{equation}
\unif{\K} \doteq\frac{\one_\K}{\vol\:\K} \label{eq:uniform}
\end{equation}
where $\one_\K$ is the indicator function of the set $\K$ defined in \eqref{indic}.
Then, the idea at the basis of the proposed method is to use a PSS approximation of the set $\K$ or, equivalently, a polynomial over approximation of the indicator function $\one_\K$,  obtained employing the framework introduced  in Sections 2 and 3.

To this end, given a degree $d \in \mathbb N$, consider the optimization problem~\eqref{l1} and let $p_d^*$ be a polynomial that achieves the optimum. If one examines  the proof of Theorem~\ref{cvg}, one can see that this polynomial has the following properties
\begin{enumerate}[i)]
\item $p_d^* \geq \one_\K$ {on} $\B$
\item As $d\rightarrow \infty$, $p_d^* \rightarrow \one_\K$ both in \ $L^1$ and almost uniformly on $\B$.
\end{enumerate}

Hence, $p_d^*$ can arbitrarily approximate (from above) the indicator function of the set $\K$., and therefore it represents a so-called ``dominating density'' of the uniform density $\unif{\K}$ {on} $\B$. More formally, there exists a value $\beta> 0$ such that $\beta p_{d}^{*}(x)\ge \unif{\K}(x)$ for all $x\in\B$.
Hence, the rejection method 
from a dominating density, discussed for instance in \cite[Section 14.3.1]{TeCaDa:13}, can be applied leading to the {following} random sampling procedure. 
\begin{algorithm}[h!]
\caption{ Uniform Sample Generation in Semialgebraic Set $\K$}
\label{alg1}
\begin{enumerate}
  \item[] {Given $d \in \mathbb N$,} let $p^*_d$ be {a} solution of 
\begin{equation}\label{l1p}
\begin{array}{ll}
\displaystyle\min_{p \in \Poly{d}} &\displaystyle \int_\B p(x)dx \\
\mathrm{s.t.} & p \geq 1 \:\:\mathrm{on}\:\: \K\\
&p \geq 0 \:\:\mathrm{on}\:\: \B.
\end{array}
\end{equation}
  \item Generate a random sample $\xi$ with density proportional to $p_{d}^{*}$ over $\B$.
  \item If $\xi\not\in\K$ go to step 1.
  \item Generate a sample $u$ uniform on $[0,\,1]$.
  \item If $u\,p_{d}^{*}(\xi)\le 1$ return $x=\xi$, else go to step 1.
\end{enumerate}
\end{algorithm}

\begin{figure}[ht!]
\centerline{
\includegraphics[width=9cm]{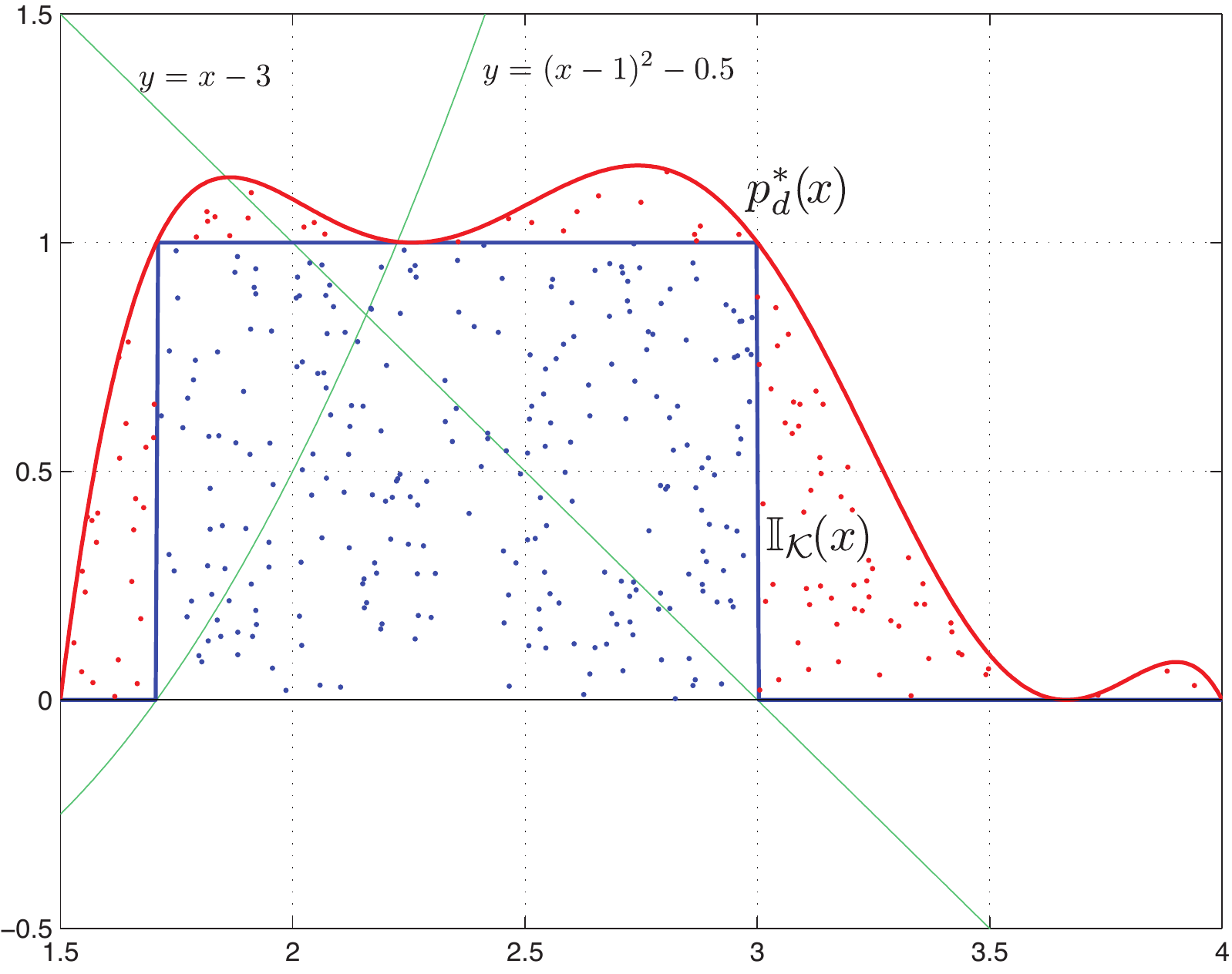}}
\caption{Illustration of the behavior of Algorithm 1 in the one-dimensional case. Blue dots are accepted samples, red dots
are rejected samples. \label{fig:oneD}}
\end{figure}

A graphical interpretation of the algorithm  is provided in Figure \ref{fig:oneD}, for the 
case of a simple one-dimensional set 
\[
\K = 
\left\{x\in\mathbb{R}\, : \,
 (x-1)^2-0.5 \ge 0,
x-3\le 0 \right\}.
\]
First, problem (\ref{l1}) is solved (for $d=8$ and $\B=[1.5,\,4]$), yielding the optimal solution
\[
p_{d}^{*}(x) =
0.069473 x^{8}
      -2.0515 x^{7}+
       23.434 x^{6}
       -139.5 x^{5}+       477.92 x^{4}
      -961.88 x^{3}+
       1090.8 x^{2}
      -606.07 x+
       107.28.
\]
As it can be seen in Figure~\ref{fig:oneD}, $p_{d}^{*}$ is ``dominating'' the indicator function $\one_{\K}$ {on} $\B$.
Then, uniform random samples are drawn in the hypograph of $p_{d}^{*}$. This is done
 by generating uniform samples $\xi$ distributed according to a probability density function (pdf)
proportional to $p_{d}^{*}$ (step 2), 
and then selecting its vertical coordinate uniformly in the interval $[0,\, \xi]$ (step 3).
Finally, if this sample falls below the indicator function $\one_{\K}$ (blue dots) it is accepted, otherwise it is rejected 
(red dots) and the process starts again.

It is intuitive that this algorithm should outperform classical rejection from the bounding set $\B$, since more importance
is given to the samples inside $\K$ through the function $p_{d}^{*}$.
To formally analyze the performance of Algorithm 1, we define the
\textit{acceptance rate} (see e.g.\ \cite{Devroye:97}) as {the reciprocal of} the  expected number
of samples that have to be drawn from $p_{d}^{*}$ in order
to find one ``good" sample, that is a sample uniformly distributed in $\K$. 
Then, the following result, which is the main theoretical result of this section, 
provides the acceptance rate of the proposed algorithm.

\smallskip

\begin{theorem}
Algorithm 1 returns a sample uniformly distributed in $\K$. Moreover, the acceptance rate 
of the algorithm is given by 
\[
\gamma_{d}=\frac{\vol\:K}{w^*_d},
\]
where $w^*_d\doteq\int_\B p^*_d(x)dx$ is the optimal solution of problem (\ref{l1}).
\end{theorem}
\vskip 3mm
\noindent
{\bf Proof:}
To prove the statement, we first note  that polynomial $p_{d}^{*}$ defines a density 
\begin{equation}
\label{poly-density}
f \doteq\frac{p_{d}^{*}}{w^*_d}
\end{equation}
over $\B$. 
Moreover, by construction, we have $p_{d}^{*}\ge\one_{\K}$ {on $\B$}, and hence 
\begin{align}
\frac{p_{d}^{*}}{w^*_d\:\vol\:\K} &\ge \frac{\one_{\K}}{w^*_d\:\vol\:\K} \label{eq:rejrate}\\
\frac{f}{\vol\:\K} &\ge  \frac{\unif{\K}f}{w^*_d}  \ge \gamma_{d} \unif{\K} \nonumber
\end{align}
{on $\B$.}
Then, it can be immediately seen that Algorithm 1 is a restatement of the classical Von Neumann rejection algorithm,
see e.g.\ \cite[Algorithm 14.2]{TeCaDa:13}, whose acceptance rate is given by the value of $\gamma_{d}$ such that 
(\ref{eq:rejrate}) holds, see for instance \cite{Devroye:86}.
\hfill$\square$\\

\smallskip
It follows that the efficiency of the random sample generation increases as $d$ increases, and becomes optimal as $d$ goes to infinity,
as reported in the next corollary.

\smallskip

\begin{corollary}
{In Algorithm 1, the acceptance rate tends to one when  increasing the degree of the polynomial approximation}, i.e.
\[
\lim_{d\to\infty} \gamma_{d} = 1.
\] 
\end{corollary}

\smallskip

Therefore, a trade-off exists   between the complexity of computing a good approximation ($d$ large) on the one hand, and 
having to wait a long time to get a ``good'' sample ($\gamma$ large), on the other hand. Note, however, that the first step can be computed off-line
for a given set $\K$, and then the corresponding polynomial $p_{d}^{*}$ can be used for efficient on-line sample generation.
Finally, we highlight that,  in order to apply Algorithm 1 in an efficient way (step 2), a computationally efficient scheme for generating random samples according to a polynomial density is required. This is discussed next.

\subsection{Sample generation from a polynomial density} \label{sec:pol_sample}

To generate a random sample according to the multivariate polynomial density $f$ defined in \eqref{poly-density}, one can {use} the so-called conditional density method described in \cite{Devroye:86}. 
This is a recursive method in which the individual entries of the multivariate samples
are generated according to their conditional probability density. We now elaborate on this. {We should note that the approach developed in this paper only provides the density up to a multiplying constant. However, to simplify the exposition to follow, we proceed as if the polynomial given is indeed a probability density function.

Assume that the bounding set is a hyperrectangle $\B=[a,b]$ of the form \eqref{Bab} 
and that we have a polynomial density $p$.  We start by computing the marginal density 
\[
p_1 : {x_1 \:\:\mapsto\:\:}  \int_{a_2}^{b_2} \cdots  \int_{a_n}^{b_n} p({x_1,x_2,\ldots,x_n})\  dx_2 \cdots dx_n
\]
and, for each $i=2,\ldots,n$ and given $\bar{x}_1, \dots \bar{x}_{i-1}$, compute conditional  marginal densities
\[
p_i : {x_i \:\:\mapsto\:\:} \int_{a_{i+1}}^{b_{i+1}} \cdots  \int_{a_n}^{b_n} p(\bar{x}_1,\ldots,\bar{x}_{i-1},x_i,x_{i+1},\ldots,x_n) \ dx_{i+1} \cdots  dx_n
\]
and respective (polynomial) cumulative distributions {$F_i$} satisfying
\[
\frac{dF_i}{dx_i} = p_i.
\]
The sampling procedure then starts by computing a sample $\bar{x}_1$ according to $F_1$ and, iteratively, computing samples $\bar{x}_{i}$  given  $\bar{x}_1, \dots \bar{x}_{i-1}$ according to the distribution $F_i$. The exact description of this procedure is described in Algorithm~\ref{alg2}. One should note that, given the density $p$, a closed form is available for all marginal and conditional densities. In other words, none of the integrations mentioned above needs to be computed numerically.
}

\begin{algorithm}[h!]
\caption{Generation from a polynomial density}
\label{alg2}
Returns a sample in $\B_{[a,b]}$ with density proportional to the polynomial 
\begin{equation}
\label{palpha2}
p : {\:\:x_1,\ldots,x_n\:\:\mapsto\:\:}
\sum_{j=1}^{n_{\alpha}}p_{j}
\prod_{\ell=1}^{n} x_{\ell}^{\alpha_{j,\ell}}
\end{equation}
\begin{enumerate}
\item
Let $i=1$ 
\item
Compute the univariate polynomial
\begin{equation}
\label{Fi}
F : {\:\:x_{i}\:\:\mapsto\:\:} \sum_{j=1}^{n_{\alpha}} \gamma_{i,j}(\bar x_{1},\ldots,\bar x_{i-1}) x_{i}^{\alpha_{j,\ell}+1}
\end{equation}
where
\begin{equation}
\label{gamma_j}
\gamma_{i,j}(\bar x_{1},\ldots,\bar x_{i-1})=
\frac{1}{a_{j,i}+1}
p_{j}\left(
\prod_{\ell=1}^{i-1} \bar x_{\ell}^{\alpha_{j,\ell}}
\right)
\left(
\prod_{\ell=j+1}^{n} 
\frac{1}{\alpha_{j,\ell}+1}
\left(
b_{\ell}^{\alpha_{j,\ell}+1}-
a_{\ell}^{\alpha_{j,\ell}+1}
\right)
\right)
\end{equation}
\item Generate a random variable $w$ uniform on $[F(a_{i}),\,F(b_{i})]$
\item Compute the unique root $\xi_{i}$ in $[a_{i},b_{i}]$ of the polynomial
${x_i\:\:\mapsto\:\:} F(x_i) -w$
\item Let $\bar x_{i}=\xi_{i}$
\item If $i<n$ let $i=i+1$ and go to (2)
\item Return $\bar x$
\end{enumerate}
\end{algorithm}

\subsection{Numerical example: sampling in a nonconvex semialgebraic set}

\begin{figure}[!ht]
\centerline{
\includegraphics[width=10cm]{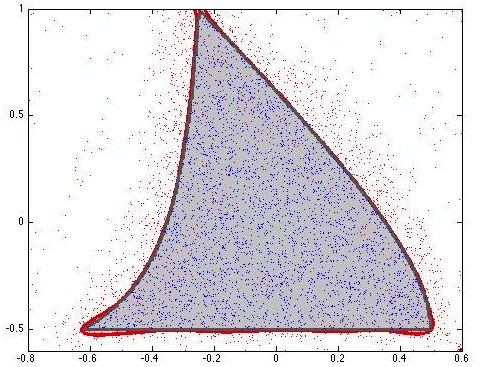}
}
\caption{Uniform random samples generated according to Algorithms 1 and 2.
The light blue area is the set $\K$ defined in \eqref{K-didier}, the pink area is the {PSS} 
$\U(p_d^*)$. The red dots are the discarded samples. The remaining samples (blue) are uniformly distributed inside $\K$.  \label{fig:didier20samples}}
\end{figure}

To demonstrate the behavior of Algorithms 1 and 2,  we revisit Example \ref{example-didier}, and generate uniform samples in the semialgebraic set $\K$ defined in \eqref{K-didier}. 
As already shown in Figure \ref{fig:schur4}, the indicator function $\one_\K$ is well approximated from above by the optimal PSS $p_{d,d}^{*}$ for $d=20$. The results of Algorithm 1 are reported in Figure \ref{fig:didier20samples}. The red points represent the points which have been discarded. To this regard, it is important to notice that also some point falling inside $\K$ has been rejected. 
This is fundamental to guarantee uniformity of the discarded points.

\section{Concluding Remarks} \label{sec:conclusion}

In this paper we have introduced the concept of polynomial superlevel sets (PSS) as a tool to construct  ``simple'' approximations of complex semialgebraic sets. Algorithms are provided for computing these approximations. Moreover, it is shown how this concept can be used to solve two important problems: i) reconstruction/approximation of sets from samples  and ii) generation of uniform samples in basic semialgebraic sets. Examples of the application of these ideas to problems in control engineering are also described. Note that the methods provided in this paper can be used to obtain probabilistic approximations of difficult sets, in the spirit of what is discussed in \cite{DaLaSh:10}. Also, in \cite{DHLS:15} the application of minimum size PSS to the approximation of the one-step reachable set of a nonlinear discrete-time function is presented, with an extension to nonlinear set filtering. Finally, we note that similar techniques can also be used to approximate transcendental (i.e. non-semi-algebraic) sets arising in systems control, e.g. regions of attraction, maximum positively invariant sets, and controllability regions.


\end{document}